\documentclass[final,leqno,onefignum]{siamltex704}

\usepackage{mathrsfs}
\usepackage{amsmath,amssymb}
\usepackage{enumitem}
\usepackage{amsfonts}
\usepackage{latexsym}
\usepackage{pictex}
\usepackage{color}
\usepackage{graphicx}
\usepackage{eucal}

\newtheorem{remark}{Remark}[section]
\newtheorem{assumption}{Assumption}[section]

\def\pseudor{pseudodifferential operator{}}

\def\rhs{right-hand side}

\def\lb{\label}

\def\g{\gamma}    
    
\def\d{\delta}

\def\p{\psi}      
   
 \def\cL{{\cal L}}   
\def\l{\lambda}   
\def\L{\Lambda}   
       
\def\n{\nu}       
      
\def\s{\sigma}

\def\Z{{\Bbb Z}}
\def\R{{\Bbb R}}

\def\qq{\quad}
\newcommand{\ma}{\begin{pmatrix}}
\newcommand{\am}{\end{pmatrix}}
\newcommand{\ca}{\begin{cases}}
\newcommand{\ac}{\end{cases}}
\let\ge\geqslant
\let\le\leqslant
\let\geq\geqslant
\let\leq\leqslant
\def\ma{\left(\begin{array}{cc}}
\def\am{\end{array}\right)}

\let\geq\geqslant
\let\leq\leqslant
\def\[{\begin{equation}}
\def\]{\end{equation}}

\def\/{\over}

\newcommand{\pdpd}[2]{\frac{\partial #1}{\partial #2}}

\newcommand{\ii}{\mathrm{i}}

\def \mr {{\mathbb R}}

\def \mn {{\mathbb N}}

\def \ha{ {\frac{1}{2}}}

\def \p {\partial}

\def \rao#1 {\frac{\p}{\p #1} #1}

\newcommand{\vphi}{\varphi}

\newcommand{\beq}{\begin{equation}}
  \newcommand{\eeq}{\end{equation}}
\newcommand{\mstrut}[1]{\mbox{\rule{0mm}{#1}}}

\def \mr {{\mathbb R}}

\def \mn {{\mathbb N}}

\def \hl{\hat{\lambda}}
\def \hm{\hat{\mu}}

\DeclareMathAlphabet{\mathcal}{OMS}{cmsy}{m}{n}

\date{\today}

\title{Semiclassical inverse spectral problem for elastic Rayleigh
  waves in isotropic media}

\author{Maarten V. de Hoop \thanks{Simons Chair in Computational and
    Applied Mathematics and Earth Science, Rice University, Houston,
    TX 77005, USA, (\texttt{mdehoop@rice.edu})}
\and Alexei Iantchenko \thanks{Department of Materials Science and
  Applied Mathematics, Faculty of Technology and Society, Malm\"{o}
  University, SE-205 06 Malm\"{o}, Sweden, (\texttt{ai@mau.se})}
\and Robert D. van der Hilst \thanks{Department of Earth, Atmospheric
  and Planetary Sciences, Massachusetts Institute of Technology,
  Cambridge, MA 02139, USA (hilst@mit.edu)}
\and Jian Zhai \thanks{Institute for Advanced Study, The Hong Kong
  University of Science and Technology, Hong Kong, China,
  (\texttt{jian.zhai@outlook.com})}}

\begin{document}

\maketitle

\pagestyle{myheadings}
\thispagestyle{plain}
\markboth{DE HOOP, IANTCHENKO, VAN DER HILST and ZHAI}{Semiclassical
  inverse spectral problem for Rayleigh waves}

\begin{abstract} 
We analyze the inverse spectral problem on the half line associated
with elastic surface waves. Here, we extend the treatment of Love
waves \cite{dHIvdHZ} to Rayleigh waves. Under certain conditions, and
assuming that the Poisson ratio is constant, we establish uniqueness
and present a reconstruction scheme for the \textit{S}-wave speed with
multiple wells from the semiclassical spectrum of these waves.
\end{abstract}

\section{Introduction}

We analyze the inverse spectral problem on the half line associated
with elastic surface waves. We discussed Love waves in a previous
paper \cite{dHIvdHZ}, and in this paper we analyze this inverse
problem for Rayleigh waves.

We study the elastic wave equation in $X = \mr^2 \times
(-\infty,0]$. In coordinates,
$$
   (x,z) ,\quad x = (x_1,x_2) \in \mr^2 ,\
                         z \in \R^{-} = (-\infty,0],
$$
we consider solutions, $u = (u_1,u_2,u_3)$, satisfying the Neumann
boundary condition at $\partial X = \{z=0\}$, to the system
\begin{equation}\label{elaswaeq}
\begin{split}
   \partial^2_t u_i + M_{il} u_l &= 0 ,\\
   u(t=0,x,z) &= 0,~~\partial_tu(t=0,x,z)=h(x,z) ,\\[-0.1cm]
   \frac{c_{i3kl}}{\rho}\partial_k u_l(t,x,z=0) &= 0 ,
\end{split}
\end{equation}
where
\begin{multline*}
M_{il} =
-\frac{\partial}{\partial z}\frac{c_{i33l}(x,z)}{\rho(x,z)}
\frac{\partial}{\partial z}     - \sum_{j,k=1}^{2}
\frac{c_{ijkl}(x,z)}{\rho(x,z)} \frac{\partial}{\partial x_j}
\frac{\partial}{\partial x_k}
- \sum_{j=1}^{2}
\frac{\partial}{\partial x_j}\frac{c_{ij3l}(x,z)}{\rho(x,z)}
\frac{\partial}{\partial z}
\\
- \sum_{k=1}^{2}
\frac{c_{i3kl}(x,z)}{\rho(x,z)} \frac{\partial}{\partial z}
\frac{\partial}{\partial x_k}
-\sum_{k=1}^{2}
\left( \frac{\partial}{\partial z} \frac{c_{i3kl}(x,z)}{\rho(x,z)} \right)
\frac{\partial}{\partial x_k}
- \sum_{j,k=1}^{2}
\left( \frac{\partial}{\partial x_j} \frac{c_{ijkl}(x,z)}{\rho(x,z)} \right)
\frac{\partial}{\partial x_k}.
\end{multline*}
Here, the stiffness tensor, $c_{ijkl}$, and density, $\rho$, are
smooth and obey the following scaling: Introducing $Z =
\frac{z}{\epsilon}$,
$$
   \frac{c_{ijkl}}{\rho}(x,z)
     = C_{ijkl}\left(x,\frac{z}{\epsilon}\right) ,
                            ~~\epsilon\in(0,\epsilon_0] ;
$$
$$
   C_{ijkl}(x,Z) = C_{ijkl}(x,Z_I) = C_{ijkl}^I(x),\quad
                           Z \leq Z_I<0 .
$$
As discussed in \cite{dHINZ}, surface waves travel along the surface
$z = 0$.

\medskip\medskip

\noindent
The remainder of the paper is organized as follows. In Section
\ref{semi-R}, we give the formulation of the inverse problems as an
inverse spectral problem on the half line and treat the simple case of
recovery of a monotonic wave-speed profile. In Section \ref{BSrule-R},
we discuss the relevant Bohr-Sommerfeld quantization, which is the
main result of this paper as it forms the key component in the study
of the inverse spectral problem. In Section \ref{inverse-R}, we give
the reconstruction scheme under appropriate assumptions, which is an
adaptation of the method of Colin de Verdi\`{e}re \cite{CdV2011}.

\section{Semiclassical description of Rayleigh waves}
\label{semi-R}

\subsection{Surface wave equation, trace and the data}

We briefly summarize the semiclassical description of elastic surface
waves \cite{dHINZ}. The leading-order symbol associated with
$M_{il}$ above is given by
\begin{multline}\label{H0}
H_{0,il}(x,\xi) =
-\frac{\partial}{\partial Z}C_{i33l}(x,Z)
\frac{\partial}{\partial Z}
\\
- \ii\sum_{j=1}^{2}
C_{ij3l}(x,Z) \xi_j
\frac{\partial}{\partial Z}
- \ii\sum_{k=1}^{2}
C_{i3kl}(x,Z) \frac{\partial }{\partial Z}
\xi_k
- \ii \sum_{k=1}^{2}
\left( \frac{\partial}{\partial Z} C_{i3kl}(x,Z) \right)
\xi_k 
\\
+ \sum_{j,k=1}^{2}
C_{ijkl}(x,Z) \xi_j \xi_k .
\hspace*{4.0cm}
\end{multline}
We view $H_0(x,\xi)$ as an ordinary differential operator in $Z$, with
domain
$$
   \mathcal{D} = \left\{ v \in H^2(\mathbb{R}^-)\ \bigg|\
      \sum_{l=1}^3\left(C_{i33l}(x,0)
    \frac{\partial v_l}{\partial Z}(0)
       + \ii \sum_{k=1}^2 C_{i3kl} \xi_k v_l(0) \right) = 0 \right\} .
$$
For an isotropic medium we have
$$
   C_{ijkl} = \hat{\lambda} \delta_{ij} \delta_{kl}
       + \hat{\mu} (\delta_{ik} \delta_{jl} + \delta_{il} \delta_{jk}) ,
$$
where $\hat{\lambda} = \frac{\lambda}{\rho},$ $\hat{\mu} =
\frac{\mu}{\rho}$, and $\lambda, \mu$ are the two Lam\'{e} moduli. The
\textit{P}-wave speed, $C_P$, is then $C_P =
\sqrt{\hat{\lambda}+2\hat{\mu}}$ and the \textit{S}-wave speed, $C_S$,
is then $C_S = \sqrt{\hat{\mu}}$. We introduce
$$
   P(\xi) = \left(\begin{array}{ccc}
        |\xi|^{-1}\xi_2& |\xi|^{-1}\xi_1 & 0 \\
        -|\xi|^{-1}\xi_1 &|\xi|^{-1}\xi_2 & 0 \\
        0 & 0 & 1 \end{array}\right) .
$$
Then 
$$
   P(\xi)^{-1}H_0(x,\xi)P(\xi) = \left(\begin{array}{cc}
   H_0^L(x,\xi) & \\ & H_0^R(x,\xi) \end{array}\right) ,
$$
where
\begin{equation} \label{Hamiltonian_Rayleigh_0}
   H_0^R(x,\xi)
   \left(\begin{array}{c} \varphi_2 \\ \varphi_3
                          \end{array}\right)
   = \left(\begin{array}{c}
-\frac{\partial}{\partial Z}(\hat{\mu}\frac{\partial\varphi_2}{\partial Z})-\ii|\xi|\left(\frac{\partial}{\partial Z}(\hat{\mu}\varphi_3)+\hat{\lambda}\frac{\partial}{\partial Z}\varphi_3\right)+(\hat{\lambda}+2\hat{\mu})|\xi|^2\varphi_2\\
-\frac{\partial}{\partial Z}\left((\hat{\lambda}+2\hat{\mu})\frac{\partial\varphi_3}{\partial Z}\right)-\ii|\xi|\left(\frac{\partial}{\partial Z}(\hat{\lambda}\varphi_2)+\hat{\mu}\frac{\partial}{\partial Z}\varphi_2\right)+\hat{\mu}|\xi|^2\varphi_3
\end{array}\right)
\end{equation} 
supplemented with Neumann boundary condition
\begin{eqnarray}
   \ii |\xi| \vphi_3(0) + \pdpd{\vphi_2}{Z}(0) &=& 0 ,
\label{Rayleighboundary1}
\\
   \ii \hat{\lambda} |\xi| \vphi_2(0)
         + (\hat{\lambda} + 2\hat{\mu}) \pdpd{\vphi_3}{Z}(0) &=& 0
\label{Rayleighboundary2}
\end{eqnarray}
for Rayleigh waves.

We assume that $H_0^R(x,\xi)$ has $\mathfrak{M}(x,\xi)$ simple
eigenvalues in its discrete spectrum
\[
   \Lambda_0< \Lambda_1 < \cdots < \Lambda_\alpha < \cdots
              < \Lambda_{\mathfrak{M}}
\]
with eigenfunctions $\Phi_{\alpha,0}(Z,x,\xi)$.
(We note the difference in labeling as
  compared with \cite{dHINZ, dHIvdHZ}.) We note, here, that
$\mathfrak{M}(x,\xi)$ increases as $|\xi|$ increases. By \cite[Theorem
  2.1]{dHINZ}, we have
\begin{equation}\label{eq:comm_sys}
   H_0^R \circ \Phi_{\alpha,0}\, 
   = \Phi_{\alpha,0}
     \circ \Lambda_\alpha+\mathcal{O}(\epsilon) .
\end{equation}
Defining
\begin{equation}\label{defchi}
   J_{\alpha,\epsilon}(Z,x,\xi) =
      \frac{1}{\sqrt{\epsilon}}\Phi_{\alpha,0}(Z,x,\xi) ,
\end{equation}
microlocally (in $x$), we can construct approximate constituent
solutions of the system \eqref{elaswaeq}, with initial values
$$
   h(x,\epsilon Z) = \sum_{\alpha=0}^{\mathfrak{M}}
   J_{\alpha,\epsilon}(Z,x,\epsilon D_x) W_{\alpha,\epsilon}(x,Z).
$$
We let $W_{\alpha,\epsilon}$ solve the initial value problems (up to
leading order)
\begin{eqnarray}
   [\epsilon^2 \p_t^2
       + \Lambda_{\alpha}(x,D_x)]
	 W_{\alpha,\epsilon}(t,x,Z) &=& 0 ,
\label{wqeaj_sys}\\
   W_{\alpha,\epsilon}(0,x,Z) &=& 0 ,\quad
   \p_t W_{\alpha,\epsilon}(0,x,Z)
                       = J_{\alpha,\epsilon} W_{\alpha}(x,Z) ,
\label{eq:initrange}
\end{eqnarray}
$\alpha = 1,\ldots,\mathfrak{M}$. We let
$\mathcal{G}_0(Z,x,t,Z',\xi;\epsilon)$ be the approximate Green's
function (microlocalized in $x$), up to leading order, for Rayleigh
waves. We may write \cite{dHINZ}
\begin{equation} \label{G0}
   \mathcal{G}_0(Z,x,t,Z',\xi;\epsilon)
   = \sum_{\alpha=0}^{\mathfrak{M}} J_{\alpha,\epsilon}(Z,x,\xi)
   \left(\frac{\mathrm{i}}{2}\mathcal{G}_{\alpha,+,0}(x,t,\xi,\epsilon)
   -\frac{\mathrm{i}}{2}\mathcal{G}_{\alpha,-,0}(x,t,\xi,\epsilon)\right)
   \Lambda^{-1/2}_\alpha(x,\xi)J_{\alpha,\epsilon}(Z',x,\xi) ,
\end{equation}
where $\mathcal{G}_{\alpha,\pm,0}$ are Green's functions for half
``wave'' equations associated with
(\ref{wqeaj_sys})-(\ref{eq:initrange}). We have the trace
\begin{equation} \label{trace_from_normal_modes}
   \int_{\mathbb{R}^-}
   \widehat{\epsilon
        \partial_t \mathcal{G}_0}(Z,x,\omega,Z,\xi;\epsilon)
        \mathrm{d} (\epsilon Z)
   = \sum_{\alpha=0}^{\mathfrak{M}}
     \delta(\omega^2 - \Lambda_\alpha(x,\xi))
        \Lambda^{1/2}_\alpha(x,\xi) + {\mathcal O}(\epsilon^{-1})
\end{equation}
from which we can extract the eigenvalues
$\Lambda_\alpha,\,\alpha=1,2,\cdots,\mathfrak{M}$ as functions of
$\xi$. We use these to recover the profile of $\hat{\mu} = C_S^2$
under

\medskip

\begin{assumption} \label{assu2}
Poisson's ratio $\nu$, with $\hl = \frac{2\nu}{1-2\nu} \hm$, of the
elastic solid is constant.
\end{assumption}

\medskip

\noindent
For a Poisson solid, $\nu = \tfrac{1}{4}$. However, we only assume
that $\nu$ is known. We may thus express $\hat{\lambda}$ in terms of
$\hat{\mu}$.

\subsection{Semiclassical spectrum}

We suppress the dependence on $x$ from now on, and introduce $h =
|\xi|^{-1}$ as another semiclassical parameter. We introduce $H_{0,h}
= h^2 H_0^R(\xi)$, that is,
\begin{equation}
   H_{0,h} \left(\begin{array}{c} \varphi_2 \\ \varphi_3
                                  \end{array}\right)
   = \left(\begin{array}{c}
   -h^2 \frac{\partial}{\partial Z}
        \left(\hat{\mu} \frac{\partial \varphi_2}{\partial Z}\right)
   - \ii h \left(\frac{\partial}{\partial Z}(\hat{\mu} \varphi_3)
   + \hat{\lambda} \frac{\partial}{\partial Z} \varphi_3\right)
   + (\hat{\lambda} + 2 \hat{\mu}) \varphi_2
\\
   -h^2 \frac{\partial}{\partial Z}
        \left((\hat{\lambda} + 2 \hat{\mu})
              \frac{\partial \varphi_3}{\partial Z}\right)
   - \ii h \left(\frac{\partial}{\partial Z}(\hat{\lambda} \varphi_2)
   + \hat{\mu} \frac{\partial}{\partial Z} \varphi_2\right)
   + \hat{\mu} \varphi_3
   \end{array}\right) ,
\lb{semop}
\end{equation}
which has eigenvalues $\l_\alpha(h) = h^2 \L_\alpha$. We invoke

\medskip

\begin{assumption} \label{assu1}
For all $Z \leq Z_I$, $\hm(Z) = \hm(Z_I)$ and $\hl(Z) =
\hl(Z_I)$. Moreover,
\begin{align}
   & 0 \ \, <\ \, \hm(0) = \, \inf_{Z \le 0} \hm(Z)
       \ \, <\ \, \hm_I = \, \sup_{Z \leq 0} \hm(Z) = \hm(Z_I) ,
\label{ass1}\\
   & \text{for all}\ Z \in [Z_I,0]\quad\text{we have}\ \
     \hat{\lambda}(Z) + 2 \hat{\mu}(Z) \geq \hat{\mu}(Z_I) .
\label{ass2}
\end{align}
\end{assumption}

\noindent
The assumption that $\hm$ attains its minimum at the boundary and its
maximum in the deep zone ($Z \le Z_I$, cf.~(\ref{ass1})) is realistic
in seismology. We write $E_0 = \hm(0)$.

\begin{remark} 
We note that if Assumption $\ref{assu2}$ is satisfied, then
$(\ref{ass2})$ requires that
\[
   2 \frac{1 - \nu}{1 - 2 \nu} \hm(Z) \ge \hm(Z_I)\qq
                               \text{for all}\,\, Z\in[Z_I,0] .
\]
\end{remark}

The spectrum of $H_{0,h}$ is divided in two parts,
\[
   \sigma(H_{0,h}) = \sigma_{\mathrm{d}}(H_{0,h})
               \cup \sigma_{\mathrm{ess}}(H_{0,h}) ,
\]
where the discrete spectrum $\sigma_{\mathrm{d}}(H_{0,h})$ consists of
a finite number of eigenvalues in $(E_0,\hm_I)$ and a lowest (subsonic) eigenvalue $\lambda_0(h)<E$, that is,
\begin{gather*}
   \lambda_0(h) < E_0 < \lambda_1(h) < \lambda_2(h) < \ldots
                 < \lambda_{\mathfrak{M}}(h) < \hat{\mu}_I ,
\end{gather*} 
and the essential spectrum $ \sigma_{\mathrm{ess}}(H_{0,h}) =
[\hm_I,\infty)$ \cite{dHINZ}. (The essential spectrum is not
  absolutely continuous for Rayleigh wave operator.)
The lowest (subsonic) eigenvalue,
  $\lambda_0(h)$, lies below $\hat{\mu}(0)$ for $h$ sufficiently
  small. Its existence and uniqueness under certain conditions (which
  are satisfied, here) are explained in \cite[Theorem~4.3]{dHINZ}. See
  also the discussion in Section~\ref{Sep_cluster}. No such phenomenon
  occurs in the case of Love waves.  Again, the number of
eigenvalues, $\mathfrak{M}$ increases as $h$ decreases.

We will study how to reconstruct the profile of $\hm$ using the
semiclassical spectrum as in \cite{CdV2011}

\medskip

\begin{definition}
For given $E$ with $E_0 < E \leq \hm(Z_I)$ and positive real number
$N$, a sequence $\mu_\alpha(h)$, $\alpha=0,1,2,\ldots$ is a
semiclassical spectrum of $H_{0,h}$ mod $o(h^N)$ in $]-\infty,E[$ if,
  for all $\lambda_\alpha(h) < E$,
$$
   \lambda_\alpha(h) = \mu_\alpha(h) + o(h^N)
$$
uniformly on every compact subset $K$ of $]-\infty,E[\,$.
\end{definition}

\medskip
 
In the remainder of the paper, we will prove

\medskip

\begin{theorem}
Under all the assumptions mentioned above and below, the function
$\hat{\mu}$ can be uniquely recovered from the semiclassical spectrum
of $H_{0,h}$ modulo $o(h^{5/2})$ below $\hat{\mu}_I$.
\end{theorem}

\subsection{Reconstruction of a monotonic profile}
\lb{S-decreasing}

In the case of a monotonic profile, the reconstruction of $\hm$ is
straightforward as it coincides with the corresponding reconstruction
in the case of Love waves.

\medskip

\begin{theorem} \label{CdV-decr}
Assume that $\hat{\mu}(Z)$ is decreasing in $[Z_I,0]$. Then the
asymptotics of the discrete spectra $\lambda_{\alpha}(h)$, $0 \leq
\alpha \leq \mathfrak{M}$ as $h \to 0$ determine the function
$\hat{\mu}(Z)$.
\end{theorem}

\medskip

This is a consequence of Weyl's law. For any $E < \hat{\mu}_I$, we
have the Weyl's law for Rayleigh waves \cite{dHINZ}:
\begin{multline*}
   \#\{ \l_\alpha(h) \leq E\} = \frac{1}{2\pi h}
   [{\rm Area}(\{ (Z,\zeta)\ :\
          (\hat{\lambda} + 2 \hat{\mu})(Z) (1 + \zeta^2) \leq E\})
\\
   + {\rm Area}(\{ (Z,\zeta)\ :\ \hat{\mu}(Z) (1 + \zeta^2) \leq E\})
          + o(1)] .
\end{multline*}
We note that the Weyl's law (in the leading order) does not depend on
boundary conditions
(\ref{Rayleighboundary1})-(\ref{Rayleighboundary2}). Due to Assumption
(\ref{ass2}), ${\rm Area}(\{ (Z,\zeta)\ :\ (\hat{\lambda} + 2
\hat{\mu})(Z) (1 + \zeta^2) \leq E\})=0$, and we get
\begin{equation} \label{countingf}
  \#\{ \lambda_\alpha(h) \leq E \}
  = \frac{1}{2\pi h} \left[
    {\rm Area}(\{ (Z,\zeta)\ :\ \hat{\mu}(Z) (1 + \zeta^2) \leq E\})
            + o(1)\right] .
\end{equation}
The procedure of reconstructing the function $\hat{\mu}$ from the
\rhs{} of (\ref{countingf}) is given in \cite[Theorem
  3.2]{dHIvdHZ}. It uses an analogue
  of Lemma 3.1 there:

\begin{lemma} \label{first}
The second eigenvalue, $\lambda_1(h)$, of $H_{0,h}$ satisfies $\lim_{h
  \rightarrow 0} \lambda_1(h) = E_0$.
\end{lemma}

In particular, similarly to Remark 4.1
  in \cite{dHIvdHZ}, under Assumption \ref{assu1}, using the Taylor
  expansion of $\hat{\mu}$ near the boundary in the Bohr-Sommerfeld
  quantization condition (\ref{BSbis}), we get that $\lambda_1(h) =
  E_0 + \mathcal{O}(h^{2/3})$. If $\hat{\mu}'(0) = 0$, then the same
  method would lead to $\lambda_1 = E_0 + \mathcal{O}(h)$.

\section{Bohr-Sommerfeld quantization}\label{BSrule-R}

For the reconstruction of the profile with (multiple) wells, we need
to establish the Bohr-Sommerfeld quantization rules for $H_{0,h}$. The
semiclassical spectrum of $H_{0,h}$ will be clustered for each well
(or half-well), due to the fact that eigenfunctions are
$\mathcal{O}(h^\infty)$ outside a well. We will establish the
quantization rules for the half-well case and the full-well case
separately.

\subsection{Half well}
\label{halfwell}

Here, we assume that the profile, $\hat{\mu}$, has a single half-well
connected to the boundary.  We follow Woodhouse and Kennett
\cite{KennettWoodhouse1978, Woodhouse1978} and rewrite $H_{0,h}
\varphi = E \varphi$ as a system of first-order ordinary differential
equations. We introduce
\[
   \psi_2 = \hat{\mu} \, (h \partial_Z (-\ii \varphi_2) + \varphi_3) ,
\qq
   \psi_3 = (\hat{\lambda}+2\hat{\mu}) \, h \partial_Z \varphi_3
                        - \hat{\lambda} \, (-\ii \varphi_2) .
\]
Then the eigenvalue problem attains the form
\[
\lb{eqnew}
   h \partial_Z
     \ma -\ii \varphi_2 \\ \varphi_3 \\ \psi_2 \\ \psi_3 \am
   = \left(\begin{array}{cccc}
     0 & -1 & \frac{1}{\hat{\mu}} & 0
     \\
     \frac{\hat{\lambda}}{\hat{\lambda}+2\hat{\mu}} & 0 & 0 &
           \frac{1}{\hat{\lambda}+2\hat{\mu}}
     \\
     \left(-E + \frac{(\hat{\lambda}+2\hat{\mu})^2
              - \hat{\lambda}^2}{\hat{\lambda}+2\hat{\mu}}\right) & 0 & 0 &
              - \frac{\hat{\lambda}}{\hat{\lambda}+2\hat{\mu}}
     \\
     0 & -E & 1 & 0
     \end{array}\right)
     \ma -\ii \varphi_2 \\ \varphi_3 \\ \psi_2 \\ \psi_3 \am
\]
supplemented with the (Neumann) boundary condition $\psi_2 = \psi_3 =
0$ at $Z = 0$. The eigenvalues of the matrix 
$$
   A_0^S := \left(\begin{array}{cccc}
   0 & -1 & \frac{1}{\hat{\mu}} & 0
   \\
   \frac{\hat{\lambda}}{\hat{\lambda}+2\hat{\mu}} & 0 & 0 & \frac{1}{\hat{\lambda}+2\hat{\mu}}
   \\
   \left(-E+\frac{(\hat{\lambda}+2\hat{\mu})^2-\hat{\lambda}^2}{\hat{\lambda}+2\hat{\mu}}\right)
   & 0 & 0 & -\frac{\hat{\lambda}}{\hat{\lambda}+2\hat{\mu}}
   \\
   0 & -E & 1 & 0 \end{array}\right)
$$
are
\[
\lb{evs}
   \pm \mathrm{i}\sqrt{\frac{E}{\hat{\lambda} + 2\hat{\mu}} - 1} ,\qq
   \pm \sqrt{1 - \frac{E}{\hat{\mu}}} .
\]
We assume existence of a single \textit{S} turning point corresponding
with a zero of $\sqrt{1 - \frac{E}{\hat{\mu}}}$ occurring at $Z =
Z_*$.

\begin{remark}
The existence of one turning point is guaranteed for any eigenvalue,
$E$, above $\hat{\mu}(0)$, while only the lowest eigenvalue falls
below $\hat{\mu}(0)$ (for $h$ sufficiently small
\textnormal{\cite{dHINZ}}). See also the discussion in
Section~\ref{Sep_cluster}. This lowest
  eigenvalue can be ignored.
\end{remark}

\medskip

Following \cite{KennettWoodhouse1978, Woodhouse1978}, we define the
matrix
\begin{align*}
   & G = G(\phi_1,\phi_2,h)
\\
   & = \left(\begin{array}{cccc}
   h^{1/6} \mathrm{Ai}'(-h^{-2/3} \phi_1) &
   h^{1/6} \mathrm{Bi}'(-h^{-2/3} \phi_1) & 0 & 0
\\ \\
   h^{-1/6} \mathrm{Ai}(-h^{-2/3} \phi_1) &
   h^{-1/6} \mathrm{Bi}(-h^{-2/3} \phi_1) & 0 & 0
\\ \\
   0 & 0 & h^{1/6} \mathrm{Ai}'(-h^{-2/3} \phi_2) &
   h^{1/6} \mathrm{Bi}'(-h^{-2/3} \phi_2)
\\ \\
   0 & 0 & h^{-1/6} \mathrm{Ai}(-h^{-2/3} \phi_2) &
   h^{-1/6} \mathrm{Bi}(-h^{-2/3} \phi_2)
   \end{array}\right) ,
\end{align*}
where $\mathrm{Ai}$ and $\mathrm{Bi}$ are Airy functions
\cite{AbramowitzStegun1965} and $\phi_1$ and $\phi_2$ are phase
functions; $G$ satisfies the equation
\[
   h \partial_Z G = Q G\quad\text{with}\quad
   Q = \left(\begin{array}{cccc}
       0 & \phi_1 \partial_Z \phi_1 & 0 & 0
       \\
       -\partial_Z \phi_1 & 0 & 0 & 0
       \\
       0 & 0 & 0 & \phi_2 \partial_Z \phi_2
       \\
       0 & 0 & -\partial_Z \phi_2 & 0
       \end{array}\right) .
\]
We search for solutions of (\ref{eqnew}) of the form
\[
\lb{F2}
   \left(\sum_{n=0}^\infty h^n Y^{(n)} \right) G(\phi_1,\phi_2, h) ,
\]
suppressing the dependencies on $E$ in the notation. Substituting
(\ref{F2}) into (\ref{eqnew}), we get from the leading order terms
\[
   A_0^S Y^{(0)} = Y^{(0)} Q .
\]
If we demand that $Y^{(0)}$ is non-singular, it follows that $A_0^S$
and $Q$ must have identical eigenvalues given in (\ref{evs}), which
implies that
$$
   \phi_1 (\partial_Z \phi_1)^2
          = \frac{E}{\hat{\lambda} + 2 \hat{\mu}} - 1 ,\qq
   \phi_2 (\partial_Z \phi_2)^2
          = \frac{E}{\hat{\mu}} - 1
$$
and, therefore,
\begin{equation}
   \phi_1(Z) =
   -\left(\mstrut{0.55cm}\right.\frac32 \int_0^Z
   \mathrm{i}\left(\mstrut{0.55cm}\right. 1 -
   \frac{E}{(\hat{\lambda} + 2 \hat{\mu})(y)}
   \left.\mstrut{0.55cm}\right)^{1/2}
   \mathrm{d}y\left.\mstrut{0.55cm}\right)^{2/3} ,
\qq
   \phi_2(Z) =
   \left(\mstrut{0.55cm}\right.\frac32 \int_{Z^*}^Z
   \left(\frac{E}{\hat{\mu}(y)} - 1\right)^{1/2} \mathrm{d}y
         \left.\mstrut{0.55cm}\right)^{2/3} ,
\end{equation}
where $Z^*$ is the unique \textit{S} turning point.

Next, we introduce explicit similarity transformations connecting
$A_0^S$ and $Q$. We introduce
$$
   \cL = \left(\begin{array}{cccc}
   0 & 1 - \frac{E}{\hat{\lambda}+2\hat{\mu}} & 0 & 0 \\
   1 & 0 & 0 & 0 \\
   0 & 0 & 0 & 1 - \frac{E}{\hat{\mu}} \\
   0 & 0 & 1 & 0 \end{array}\right) .
$$
Then
\[
\lb{simtr}
   R^{-1} A_0^S R = \Phi Q \Phi^{-1} = \cL ,
\]
where the similarity transformations, $R$ and $\Phi$, defined by
(\ref{simtr}) (formula (56) in \cite{Woodhouse1978}) are given by
$$
   R = \left(\begin{array}{cccc}
   1 & 0 & 0 & 1 - \frac{E}{\hat{\mu}}
   \\ \\
   0 & 1 - \frac{E}{\hat{\lambda}+2\hat{\mu}} & 1 & 0 
   \\ \\
   0 & 2\hat{\mu} \left(1 - \frac{E}{
           \hat{\lambda}+2\hat{\mu}}\right) & 2 \hat{\mu} - E & 0
   \\ \\
   2\hat{\mu} - E & 0 & 0 & 2 (\hat{\mu} - E) \end{array}\right)
$$
and
$$
   \Phi= \left(\begin{array}{cccc}|\partial_Z\phi_1|^{1/2} & 0 & 0 & 0
\\
   0 & -|\partial_Z\phi_1|^{-1/2} & 0 & 0
\\
   0 & 0 & |\partial_Z\phi_2|^{1/2} & 0
\\
   0 & 0 & 0 &-|\partial_Z\phi_2|^{-1/2}\end{array}\right) .
$$
Writing $Y^{(n)} = R T^{(n)} \Phi$, expansion (\ref{F2}) takes the
form
\begin{equation} \label{F2bis}
   R \left(\sum_{n=0}^\infty h^nT^{(n)}\right)
                   \Phi G(\phi_1,\phi_2,h) .
\end{equation}
Denoting $T = \sum_{n=1}^\infty h^{n-1}T^{(n)}$, expansion
(\ref{F2bis}) takes the form
$$
   R (I + h T) \ma E_1 & 0 \\ 0 & E_2 \am ,
$$
where
$$
   E_1=  \left(\begin{array}{cc}h^{1/6} |\partial_Z\phi_1|^{1/2} \mathrm{Ai}'(-h^{-2/3}\phi_1) &h^{1/6} |\partial_Z\phi_1|^{1/2} \mathrm{Bi}'(-h^{-2/3}\phi_1) \\ \\
-h^{-1/6}|\partial_Z\phi_1|^{-1/2}\mathrm{Ai}(-h^{-2/3}\phi_1) &-h^{-1/6}|\partial_Z\phi_1|^{-1/2}\mathrm{Bi}(-h^{-2/3}\phi_1) \end{array}\right)=:\ma a_1 &b_1\\ c_1 &d_1\am
$$
and
$$
   E_2 = \left(\begin{array}{cc}
h^{1/6} |\partial_Z\phi_2|^{1/2} \mathrm{Ai}'(-h^{-2/3}\phi_2) &h^{1/6} |\partial_Z\phi_2|^{1/2} \mathrm{Bi}'(-h^{-2/3}\phi_2)
   \\ \\
-h^{-1/6} |\partial_Z\phi_2|^{-1/2}\mathrm{Ai}(-h^{-2/3}\phi_2) &-h^{-1/6}|\partial_Z\phi_2|^{-1/2}\mathrm{Bi}(-h^{-2/3}\phi_2) \end{array}\right)
   =: \ma a_2 & b_2 \\ c_2 & d_2 \am .
$$

The matrix $R$ corresponds to a local decomposition of the
displacement field into standing \textit{P}- and \textit{S}-wave
constituents. The interactions of these standing waves with one
another and with the velocity gradient are of lower order in $h$ and
appear through the matrix $T$ (given in \cite{KennettWoodhouse1978}
for the spherical case). We note the asymptotic behavior,
$$
   h^{-1/6} |\partial_Z \phi_2(Z)|^{-1/2}
      \mathrm{Ai}(-h^{-2/3} \phi_2(Z)) \sim
   \left(\mstrut{0.55cm}\right.
   \frac{E}{\hat{\mu}(Z)} - 1\left.\mstrut{0.55cm}\right)^{-1/4}
   \cos\left(\mstrut{0.55cm}\right.\frac{1}{h} \int_{Z^*}^Z
   \left(\frac{E}{\hat{\mu}(y)} -1\right)^{1/2} \mathrm{d}y
          - \frac{\pi}{4}\left.\mstrut{0.55cm}\right) ,
$$
$$
   h^{1/6} |\partial_Z \phi_2(Z)|^{1/2}
   \mathrm{Ai}'(-h^{-2/3} \phi_2(Z)) \sim
   -\left(\mstrut{0.55cm}\right.\frac{E}{\hat{\mu}(Z)} - 1\left.\mstrut{0.55cm}\right)^{1/4}
   \sin\left(\mstrut{0.55cm}\right.\frac{1}{h} \int_{Z^*}^Z \left(\frac{E}{\hat{\mu}(y)} -1\right)^{1/2} \mathrm{d}y - \frac{\pi}{4}\left.\mstrut{0.55cm}\right) 
$$
in the allowed (propagating) region for $S$ waves ($\mathrm{Bi}$
similar), and
$$
   h^{-1/6} |\partial_Z \phi_1(Z)|^{-1/2}
   \mathrm{Ai}(-h^{-2/3} \phi_1(Z)) \sim
   \frac12\left(1-\frac{E}{(\hat{\lambda}+2\hat{\mu})(Z)}\right)^{-1/4}
   \exp\left(\mstrut{0.55cm}\right.\!-\frac{1}{h}\int_{0}^Z
   \left(\mstrut{0.55cm}\right.\!1-\frac{E}{(\hat{\lambda}+2\hat{\mu})(y)}\!\left.\mstrut{0.55cm}\right)^{1/2} \mathrm{d}y\!\left.\mstrut{0.55cm}\right)\!,
$$
$$
   h^{1/6} |\partial_Z \phi_1(Z)|^{1/2}
   \mathrm{Ai}'(-h^{-2/3} \phi_1(Z)) \sim
   -\frac12 \left(1-\frac{E}{(\hat{\lambda}+2\hat{\mu})(Z)}\right)^{1/4}
   \exp\left(\mstrut{0.55cm}\right.\!-\frac{1}{h}\int_{0}^Z
   \left(\mstrut{0.55cm}\right.\!1-\frac{E}{(\hat{\lambda}+2\hat{\mu})(y)}\!\left.\mstrut{0.55cm}\right)^{1/2} \mathrm{d}y\!\left.\mstrut{0.55cm}\right) 
$$
in the forbidden (evanescent) region for $P$ waves ($\mathrm{Bi}$
similar but exponentially increasing so that any $\mathrm{Bi}$ term
must be excluded in this region, see \cite{KennettWoodhouse1978}).

The solution is then given by (see also (11) in
\cite{KennettWoodhouse1978})
$$
   \ma -\ii\varphi_2 \\ \varphi_3 \\ \psi_2 \\ \psi_3 \am
   \sim R \, (I + h T)
   \ma E_1 & 0 \\ 0 & E_2 \am \ma c_\alpha \\ 0 \\ d_\beta \\ 0 \am .
$$ 
We calculate the zeroth order explicitly,
$$
   \ma-\ii\varphi_2\\ \\ \varphi_3\\ \\ \psi_2\\ \\ \psi_3\am\sim \left(\begin{array}{cccc}
a_1&b_1&\left(1-\frac{E}{\hat{\mu}}\right)c_2&\left(1-\frac{E}{\hat{\mu}}\right)d_2\\ \\
\left(1-\frac{E}{\hat{\lambda}+2\hat{\mu}}\right)c_1&\left(1-\frac{E}{\hat{\lambda}+2\hat{\mu}}\right)d_1&a_2&b_2\\ \\
2\hat{\mu}\left(1-\frac{E}{\hat{\lambda}+2\hat{\mu}}\right)c_1&2\hat{\mu}\left(1-\frac{E}{\hat{\lambda}+2\hat{\mu}}\right)d_1&\left(2\hat{\mu}-E\right)a_2&\left(2\hat{\mu}-E\right)b_2\\ \\
\left(2\hat{\mu}-E\right)a_1&\left(2\hat{\mu}-E\right)b_1&2(\hat{\mu}-E)c_2&2(\hat{\mu}-E)d_2\end{array}\right)\ma c_\alpha\\ \\ 0\\ \\ d_\beta\\ \\ 0\am .
$$
We get
\begin{align*}&\ma-\ii\varphi_2\\ \varphi_3\am\sim\ma a_1 & -\left(1-\frac{E}{\hat{\mu}}\right)c_2\\ \left(1-\frac{E}{\hat{\lambda}+2\hat{\mu}}\right)c_1 &a_2\am 
   \ma c_\alpha \\ d_\beta \am
\\[0.25cm]
&\!=\!\ma h^{1/6} |\partial_Z\phi_1|^{1/2} \mathrm{Ai}'(-h^{-2/3}\phi_1) & -\left(1-\frac{E}{\hat{\mu}}\right)h^{-1/6}|\partial_Z\phi_2|^{-1/2}\mathrm{Ai}(-h^{-2/3}\phi_2) \\ \\ -\left(1-\frac{E}{\hat{\lambda}+2\hat{\mu}}\right)h^{-1/6}|\partial_Z\phi_1|^{-1/2}\mathrm{Ai}(-h^{-2/3}\phi_1) &h^{1/6} |\partial_Z\phi_2|^{1/2} \mathrm{Ai}'(-h^{-2/3}\phi_2)\am
   \ma c_\alpha \\ d_\beta \am
\end{align*}
and
\begin{align*}
   &\ma \psi_2 \\ \psi_3 \am
   \sim \ma 2\hat{\mu}\left(1-\frac{E}{\hat{\lambda}+2\hat{\mu}}\right)c_1 & \left(2\hat{\mu}-E\right)a_2
\\
\left(2\hat{\mu}-E\right)a_1 &2(\hat{\mu}-E)c_2\am  \ma c_\alpha \\d_\beta \am
\\[0.25cm]
&\!=\!\ma -2\hat{\mu}\left(1-\frac{E}{\hat{\lambda}+2\hat{\mu}}\right)h^{-1/6}|\partial_Z\phi_1|^{-1/2}\mathrm{Ai}(-h^{-2/3}\phi_1) & \left(2\hat{\mu}-E\right)h^{1/6} |\partial_Z\phi_2|^{1/2} \mathrm{Ai}'(-h^{-2/3}\phi_2)\\ \\
\left(2\hat{\mu}-E\right)h^{1/6} |\partial_Z\phi_1|^{1/2} \mathrm{Ai}'(-h^{-2/3}\phi_1) &-2(\hat{\mu}-E)h^{-1/6}|\partial_Z\phi_2|^{-1/2} \mathrm{Ai}(-h^{-2/3}\phi_2) \am  \ma c_\alpha \\ d_\beta \am\!.
\end{align*}
Using the asymptotics of the Airy functions in the allowed region for
\textit{S} and in the forbidden region for \textit{P}, and imposing on the
expansion the boundary condition, $\psi_1 = \psi_3 = 0$ at $Z = 0$, we
get from the zeroth order terms in $h$,
\begin{eqnarray*}
&\begin{pmatrix}
   -\hat{\mu} \left(1 - \frac{E}{\hat{\lambda}+2\hat{\mu}}\right)
     \left(1-\frac{E}{\hat{\lambda}+2\hat{\mu}} \right)^{-1/4} &
   -(2 \hat{\mu} - E)
     \left(\frac{E}{\hat{\mu}} - 1\right)^{1/4}
           \sin\left(
    \displaystyle{\frac{1}{h}
    \int_{Z^*}^0} \left(\frac{E}{\hat{\mu}(Z)} - 1\right)^\frac12
                  \mathrm{d}Z - \frac{\pi}{4}\right) \\
   -\frac12 (2 \hat{\mu} - E)
     \left(1-\frac{E}{\hat{\lambda}+2\hat{\mu}} \right)^{1/4} &
   -2 (\hat{\mu} - E)
    \left(\frac{E}{\hat{\mu}} - 1\right)^{-1/4}
          \cos\left(
    \displaystyle{\frac{1}{h}
    \int_{Z^*}^0} \left(\frac{E}{\hat{\mu}(Z)} - 1\right)^\frac12
                 \mathrm{d}Z - \frac{\pi}{4}\right)    
   \end{pmatrix}
   \ma c_\alpha \\ d_\beta \am\\
&\qquad\qquad
   = \ma 0 \\ 0 \am .
\end{eqnarray*}

There is a non-trivial solution if
\begin{equation} \label{BS}
   \tan\left(\frac{1}{h}
   \int_{Z^*}^0 \left(\frac{E}{\hat{\mu}(Z)} - 1\right)^{1/2}
                \mathrm{d}Z - \frac{\pi}{4}\right)
   = - \frac{4 \hat{\mu}(0)(\hat{\mu}(0) - E)
      \left(1 - \frac{E}{(\hat{\lambda}+2\hat{\mu})(0)}\right)^{1/2}}{
      (2 \hat{\mu}(0) - E)^2
               \left(\frac{E}{\hat{\mu}(0)} - 1\right)^{1/2}} ,
\end{equation}
which is the implicit Bohr-Sommerfeld quantization in leading order in
$h$, sufficient for the further analysis. We note that in the allowed
region for \textit{S} and in the forbidden region for \textit{P}, the
right-hand side of (\ref{BS}) is negative. Then (\ref{BS}) implies the
Bohr-Sommerfeld quantization condition in leading order in $h$,
\begin{equation} \label{BSbis}
   \frac{1}{h} \int_{Z^*}^0 \left(\frac{E}{\hat{\mu}(Z)}
     - 1 \right)^{1/2} \mathrm{d}Z + \frac{3\pi}{4}
   + \arctan\left(\frac{4 \hat{\mu}(0)(\hat{\mu}(0) - E)
      \left(1 - \frac{E}{(\hat{\lambda}+2\hat{\mu})(0)}\right)^{1/2}}{
      (2 \hat{\mu}(0) - E)^2
      \left(\frac{E}{\hat{\mu}(0)} - 1\right)^{1/2}}\right)
   = \alpha \pi + {\mathcal O}(h) ,
\end{equation}
for $\alpha=1,2,\cdots$. The estimate ${\mathcal O}(h)$ follows from
Poincar\'{e}-type expansions of the Airy functions~\footnote{For
  Poincar\'{e}-type expansions of the Airy functions, see, for
  example, \texttt{https://dlmf.nist.gov/9.7}.}.

\subsection{Wells separated from the boundary}

\subsubsection{Diagonalization of the Rayleigh matrix operator} 

For the semiclassical wells separated from the boundary, $Z=0$, we may
apply techniques used for semiclassical matrix-valued spectral problems
on the whole line, namely semiclassical diagonalization.

The Weyl symbol of $H_{0,h}$ is given by $$\s^{\rm
  W}(H_{0,h})=q=q_0+hq_1+h^2q_2,$$ with
\begin{equation} \label{princs}
   q_0 = \left(\begin{array}{lr}  \hat{\mu}\zeta^2+(\hat{\lambda}+2\hat{\mu})&(\hat{\lambda}+\hat{\mu})\zeta\\
(\hat{\lambda}+\hat{\mu})\zeta&(\hat{\lambda}+2\hat{\mu})\zeta^2+\hat{\mu}\end{array}\right),
\
   q_1=\frac{1}{2\ii}\left(\begin{array}{cc} 0&\hat{\mu}'-\hat{\lambda}'\\
\hat{\lambda}'-\hat{\mu}'&0\end{array}\right) ,
\
   q_2 = \frac14 \left(\begin{array}{cc} \hat{\mu}''&0\\
0&\hat{\lambda}''+2\hat{\mu}''\end{array}\right)
\end{equation}
(cf.~(\ref{semop})). To prove this fact, we use the Moyal product
defined as follows (see \cite{CdV2005})
$$
a \star b:= \sum_{j=0}^\infty\frac{1}{j!}\left(\frac{h}{2\ii}\right)^j\{a,b\}_j ,$$ with
\begin{multline*}
\{a,b\}_j(Z,\zeta) :=[(\partial_\zeta\partial_{Z_1}-\partial_Z\partial_{\zeta_1})^j
a(Z,\zeta) b(Z_1,\zeta_1)]_{(Z_1,\zeta_1)=(Z,\zeta)}
\\
= \sum_{n=0}^j\ma j\\n\am (-1)^n(\partial_Z^n\partial_\zeta^{j-n}a)(Z,\zeta) \,
(\partial_Z^{j-n}\partial_\zeta^n b)(Z,\zeta)
\end{multline*}
with the property
\begin{multline*}
   a\star b=\{a,b\}_0+\frac{h}{2\ii}\{a,b\}_1
          -\frac{h^2}{8}\{a,b\}_2+{\mathcal O}(h^3) ,
\\
   \{a,b\}_0 = a b ,
\qq
   \{a,b\}_1 = \partial_\zeta a \partial_Z b - \partial_Z a \partial_\zeta b ,
\qq
   \{a,b\}_2 = \partial_\zeta^2 a \partial_Z^2 b
        - 2 (\partial_Z \partial_\zeta a)(\partial_Z \partial_\zeta b)
                      + \partial_Z^2 a \partial_\zeta^2 b .
\end{multline*}
The expressions in (\ref{princs}) follow from the calculations
below
\begin{align*}
  \zeta\star(\hat{\l}+2\hat{\mu})&=\zeta(\hat{\l}+2\hat{\mu})+\frac{h}{2\ii}(\hat{\l}'+2\hat{\mu}'),\\ \zeta\star(\hat{\l}+2\hat{\mu})\star\zeta&=\zeta(\hat{\l}+2\hat{\mu})\star\zeta+\frac{h}{2\ii}(\hat{\l}'+2\hat{\mu}')\star\zeta\\ &=\zeta^2(\hat{\l}+2\hat{\mu})-\frac{h}{2\ii}\zeta(\hat{\l}'+2\hat{\mu}')+\frac{h}{2\ii}(\hat{\l}'+2\hat{\mu}')\zeta-\left(\frac{h}{2i}\right)^2
  (\hat{\l}''+2\hat{\mu}'')\\ &=\zeta^2(\hat{\l}+2\hat{\mu})+\frac{h^2}{4}
  (\hat{\l}''+2\hat{\mu}''),\\ \hat{\mu}\star\zeta&=\hat{\mu}\zeta-\frac{h}{2\ii}\hat{\mu}'.\end{align*}

We use the method developed by Taylor \cite[Section
  3.1]{HelfferSjostrand1990} to diagonalize the matrix-valued operator
$H_{0,h}$ to any order in $h$.

\medskip

\begin{theorem}[Diagonalization]\label{Th-decoupling} There exists a unitary \pseudor{} $U$ and diagonal operator 
\begin{equation} \label{diagonal}
   \widehat{H}_{0,h} = \left(\begin{array}{cc} H_{0,h,1} & 0 \\
                             0 & H_{0,h,2} \end{array}\right) 
\end{equation}
such that
\[\lb{decoup}U^*H_{0,h} U=\widehat{H}_{0,h}+{\mathcal O}(h^\infty).\]  Here, $H_{0,h,i},$ $i=1,2,$ are \pseudor{}s with symbols  
\begin{align*}&\s^{\rm W}(H_{0,h,1}) =(\hat{\lambda}+2\hat{\mu})(1+\zeta^2)+h^2\alpha_2+\cdots,\\&\s^{\rm W}(H_{0,h,2})=\hat{\mu}(1+\zeta^2)+h^2\delta_2+\cdots \end{align*}
where
\begin{align*}\alpha_2&=\frac{\hat{\lambda}''+2\hat{\mu}''}{4}+ \frac{1}{\zeta^2+1}\left\{-\frac{1}{2}\hat{\l}''+\hat{\mu}''+\frac{4(\hat{\mu}')^2(2\hl+3\hat{\mu})}{(\hat{\lambda}+\hat{\mu})^2} \right\},\\\delta_2&=\frac{\hm''}{4}+ \frac{1}{\zeta^2+1}\left\{\frac32 \hat{\mu}''-\frac{4\hl(\hat{\mu}')^2}{(\hat{\lambda}+\hat{\mu})^2} \right\}.\end{align*}
Note that the $h^1$-order terms vanish.
\end{theorem}

\medskip

\begin{proof}
We introduce a unitary operator $U_0$, which is the Weyl quantization
of the matrix-symbol
$$
   Q  =\frac{1}{\sqrt{\zeta^2 + 1}} \left(\begin{array}{cc}
     1 & -\zeta \\ \zeta & 1 \end{array}\right) ,
\quad
   Q^{-1} = \frac{1}{\sqrt{\zeta^2 + 1}} \left(\begin{array}{cc}
     1 & \zeta \\ -\zeta & 1 \end{array}\right) ,
$$
which diagonalizes the principal symbol $q_0$, that is,
$$
   Q^{-1}q_0(Z,\zeta)Q =
   \left(\begin{array}{cc}
   (\hat{\lambda} + 2 \hat{\mu}) (1 + \zeta^2) & 0 \\
   0 & \hat{\mu} (1 + \zeta^2) \end{array}\right) .
$$
First, we calculate the $h^1$-order correction, that is the second
term in the \rhs{} of
\[\lb{abgd}
   Q^{-1}\star q\star Q = \hat{q}
   = \left(\begin{array}{cc} (\hat{\lambda} + 2 \hat{\mu}) (1+\zeta^2) & 0
     \\ 0 & \hat{\mu} (1 + \zeta^2) \end{array}\right)
   + h \left(\begin{array}{cc} \alpha_1 & \beta_1 \\ \gamma_1 & \delta_1
             \end{array}\right) + {\mathcal O}(h^2) .
\]
It will follow that $\alpha_1 = \delta_1 = 0$. Later, we will also
need the explicit form of diagonal entries of the next order
correction. Therefore, we keep the $h^2-$order terms in our
calculations. We introduce
\begin{align*}
\kappa_1 & := \left(\frac{1}{\sqrt{\zeta^2+1}}\right)'=-\frac{\zeta}{(\zeta^2+1)\sqrt{\zeta^2+1}}, \\
\kappa_2 & := \left(\frac{\zeta}{\sqrt{\zeta^2+1}}\right)'=\frac{1}{(\zeta^2+1)\sqrt{\zeta^2+1}}, \\
\kappa_3 & := \left(\frac{1}{\sqrt{\zeta^2+1}}\right)''=\frac{2\zeta^2-1}{(\zeta^2+1)^2\sqrt{\zeta^2+1}}, \\
\kappa_4 & := \left(\frac{\zeta}{\sqrt{\zeta^2+1}}\right)''=-\frac{3\zeta}{(\zeta^2+1)^2\sqrt{\zeta^2+1}}.
\end{align*}
 We start with the calculation of $p_0\star Q$ modulo terms of order $h^3$,
\begin{align*}
&\begin{pmatrix}
(\hat{\l}+2\hat{\mu})\sqrt{\zeta^2+1}&-\hat{\mu}\zeta\sqrt{\zeta^2+1}\\[1em]
(\hat{\l}+2\hat{\mu})\zeta\sqrt{\zeta^2+1}
&\hat{\mu}\sqrt{\zeta^2+1}
\end{pmatrix}
\\[1em]
&\qquad
   +\frac{h}{2\ii}\begin{pmatrix}
 -\left(\hat{\mu}'\zeta^2+(\hat{\l}'+2\hat{\mu}')\right)\kappa_1-(\hat{\l}'+\hat{\mu}')\zeta\kappa_2 &
\left(\hat{\mu}'\zeta^2+(\hat{\l}'+2\hat{\mu}')\right)\kappa_2-(\hat{\l}'+\hat{\mu}')\zeta\kappa_1 \\[1em] 
-\left((\hat{\l}'+2\hat{\mu}')\zeta^2+\hat{\mu}'\right)\kappa_2-(\hat{\l}'+\hat{\mu}')\zeta\kappa_1 & 
-\left((\hat{\l}'+2\hat{\mu}')\zeta^2+\hat{\mu}'\right)\kappa_1+(\hat{\l}'+\hat{\mu}')\zeta\kappa_2
\end{pmatrix}
\\[1em]
&\qquad
   -\frac{h^2}{8}\begin{pmatrix}
\left(\hat{\mu}''\zeta^2+(\hat{\l}''+2\hat{\mu}'')\right)\kappa_3+(\hat{\l}''+\hat{\mu}'')\zeta\kappa_4
&-\left(\hat{\mu}''\zeta^2+(\hat{\l}''+2\hat{\mu}'')\right)\kappa_4+(\hat{\l}''+\hat{\mu}'')\zeta\kappa_3 \\[1em] 
(\hat{\l}''+\hat{\mu}'')\zeta\kappa_3+\left((\hat{\l}''+2\hat{\mu}'')\zeta^2+\hat{\mu}''\right)\kappa_4
&-(\hat{\l}''+\hat{\mu}'')\zeta\kappa_4+\left((\hat{\l}''+2\hat{\mu}'')\zeta^2+\hat{\mu}''\right)\kappa_3
\end{pmatrix},
\end{align*}
where the second  term simplifies to
\begin{align*}
&\frac{h}{2\ii}
\begin{pmatrix} -(\hat{\l}'+2\hat{\mu}')(\zeta^2+1)\kappa_1-(\hat{\l}'+\hat{\mu}')\dfrac{\zeta}{\sqrt{\zeta^2+1}}&[\hat{\mu}'\zeta^2+(\hat{\l}'+2\hat{\mu}')]\dfrac{1}{\sqrt{\zeta^2+1}}+\hat{\mu}'\zeta(\zeta^2+1)\kappa_1\\[1em]
 -[(\hat{\l}'+2\hat{\mu}')\zeta^2+\hat{\mu}']\dfrac{1}{\sqrt{\zeta^2+1}}-(\hat{\l}'+2\hat{\mu}')\zeta(\zeta^2+1)\kappa_1& -\hat{\mu}'(\zeta^2+1)\kappa_1+(\hat{\l}'+\hat{\mu}')\dfrac{\zeta}{\sqrt{\zeta^2+1}}
\end{pmatrix}
\\
=&\frac{h}{2\ii}\left(\begin{array}{lr} \hat{\mu}'\dfrac{\zeta}{\sqrt{\zeta^2+1}}&(\hat{\l}'+2\hat{\mu}')\dfrac{1}{\sqrt{\zeta^2+1}}\\[1em]
 -\hat{\mu}'\dfrac{1}{\sqrt{\zeta^2+1}}&(\hat{\l}'+2\hat{\mu}')\dfrac{\zeta}{\sqrt{\zeta^2+1}}\end{array}\right) ,
\end{align*}
and the third term simplifies as follows,
\begin{align*}
&-\frac{h^2}{8}\left(\begin{array}{cc} [\hat{\mu}''\zeta^2+(\hat{\l}''+2\hat{\mu}'')]\kappa_3+(\hat{\l}''+\hat{\mu}'')\zeta\kappa_4
&-[\hat{\mu}''\zeta^2+(\hat{\l}''+2\hat{\mu}'')]\kappa_4+(\hat{\l}''+\hat{\mu}'')\zeta\kappa_3\\[1em]
 (\hat{\l}''+\hat{\mu}'')\zeta\kappa_3+[(\hat{\l}''+2\hat{\mu}'')\zeta^2+\hat{\mu}'']\kappa_4
&-(\hat{\l}''+\hat{\mu}'')\zeta\kappa_4+[(\hat{\l}''+2\hat{\mu}'')\zeta^2+\hat{\mu}'']\kappa_3\end{array}\right)\nonumber\\ \nonumber\\
=&-\frac{h^2}{8}\begin{pmatrix}
\dfrac{\zeta^2(2\hat{\mu}''\zeta^2-\hat{\l}'')-(\hat{\l}''+2\hat{\mu}'')}{(\zeta^2+1)^2\sqrt{\zeta^2+1}}
&\dfrac{\zeta(2\hat{\l}''+5\hat{\mu}'')}{(\zeta^2+1)\sqrt{\zeta^2+1}} \\[2em]
-\dfrac{\zeta(\hat{\l}''+4\hat{\mu}'')}{(\zeta^2+1)\sqrt{\zeta^2+1}}
&\dfrac{2\zeta^4(\hat{\l}''+2\hat{\mu}'')+\zeta^2(2\hat{\l}''+3\hat{\mu}'')-\hat{\mu}''}{(\zeta^2+1)^2\sqrt{\zeta^2+1}}
\end{pmatrix}.
\end{align*}
Thus we get for $q_0 \star Q$ modulo terms of order $h^3$,
\begin{align*}
&\begin{pmatrix}
(\hat{\l}+2\hat{\mu})\sqrt{\zeta^2+1}&-\hat{\mu}\zeta\sqrt{\zeta^2+1}\\[1em]
(\hat{\l}+2\hat{\mu})\zeta\sqrt{\zeta^2+1}&\hat{\mu}\sqrt{\zeta^2+1}
\end{pmatrix}
+\frac{h}{2\ii}
\begin{pmatrix}
\hat{\mu}'\dfrac{\zeta}{\sqrt{\zeta^2+1}}&(\hat{\l}'+2\hat{\mu}')\dfrac{1}{\sqrt{\zeta^2+1}}\\
-\hat{\mu}'\dfrac{1}{\sqrt{\zeta^2+1}}&(\hat{\l}'+2\hat{\mu}')\dfrac{\zeta}{\sqrt{\zeta^2+1}}
\end{pmatrix}\\[1em]
&\hspace*{3.5cm} -\frac{h^2}{8}
\begin{pmatrix}
\dfrac{\zeta^2(2\hat{\mu}''\zeta^2-\hat{\l}'')-(\hat{\l}''+2\hat{\mu}'')}{(\zeta^2+1)^2\sqrt{\zeta^2+1}}&
\dfrac{\zeta(2\hat{\l}''+5\hat{\mu}'')}{(\zeta^2+1)\sqrt{\zeta^2+1}}\\[2em]
 -\dfrac{\zeta(\hat{\l}''+4\hat{\mu}'')}{(\zeta^2+1)\sqrt{\zeta^2+1}}&
\dfrac{2\zeta^4(\hat{\l}''+2\hat{\mu}'')+\zeta^2(2\hat{\l}''+3\hat{\mu}'')-\hat{\mu}''}{(\zeta^2+1)^2\sqrt{\zeta^2+1}}
\end{pmatrix} .
\end{align*}
Now, we calculate $Q^{-1}\star q_0\star Q$ modulo terms of order $h^2$,
\begin{align*}
&\begin{pmatrix}
(\hat{\l}+2\hat{\mu})(\zeta^2+1)&0\\
0&\hat{\mu}(\zeta^2+1)
\end{pmatrix} 
+\frac{h}{2\ii}
\left(\begin{array}{cc} 
0&\hat{\l}'+2\hat{\mu}'\\
-\hat{\mu}'&0\end{array}\right)
\\[1em]
&\qquad +\frac{h}{2\ii}
\begin{pmatrix} 
\kappa_1(\hat{\l}'+2\hat{\mu}')\sqrt{\zeta^2+1}+\kappa_2(\hat{\l}'+2\hat{\mu}')\zeta\sqrt{\zeta^2+1}& 
\kappa_1(-\hat{\mu}')\zeta\sqrt{\zeta^2+1}+\kappa_2\hat{\mu}'\sqrt{\zeta^2+1}\\[1em] 
-\kappa_2(\hat{\l}'+2\hat{\mu}')\sqrt{\zeta^2+1}+\kappa_1(\hat{\l}'+2\hat{\mu}')\zeta\sqrt{\zeta^2+1}& 
\kappa_2\hat{\mu}'\zeta\sqrt{\zeta^2+1}+\kappa_1\hat{\mu}'\sqrt{\zeta^2+1}
\end{pmatrix} \\[1em]
=&\begin{pmatrix}
(\hat{\l}+2\hat{\mu})(\zeta^2+1)&0\\
0&\hat{\mu}(\zeta^2+1)
\end{pmatrix}+
\frac{h}{2\ii}
\begin{pmatrix} 0&\hat{\l}'+3\hat{\mu}'\\
-(\hat{\l}'+3\hat{\mu}')&0
\end{pmatrix},
\end{align*}
which together with (\ref{princs}) shows that $\alpha_1 = \delta_1 =
0$.

Then we calculate the terms of order $h^2$. There are three terms.

\medskip 

\textit{First term}: the term of order $h^2$ in 
\begin{align*}
\begin{pmatrix}
\dfrac{1}{\sqrt{\zeta^2+1}} & \dfrac{\zeta}{\sqrt{\zeta^2+1}}\\
-\dfrac{\zeta}{\sqrt{\zeta^2+1}} & \dfrac{1}{\sqrt{\zeta^2+1}}
\end{pmatrix} \star 
\begin{pmatrix}
(\hat{\l}+2\hat{\mu})\sqrt{\zeta^2+1}&-\hat{\mu}\zeta\sqrt{\zeta^2+1}\\[2em]
(\hat{\l}+2\hat{\mu})\zeta\sqrt{\zeta^2+1}&\hat{\mu}\sqrt{\zeta^2+1}
\end{pmatrix}
\end{align*}
is 
\begin{align*}
&-\frac{h^2}{8}
\begin{pmatrix}
\kappa_3(\hat{\l}''+2\hat{\mu}'')\sqrt{\zeta^2+1}+\kappa_4 (\hat{\l}''+2\hat{\mu}'')\zeta\sqrt{\zeta^2+1}&
-\kappa_3\hat{\mu}''\zeta\sqrt{\zeta^2+1}+\kappa_4 \hat{\mu}''\sqrt{\zeta^2+1}\\[1em] 
-\kappa_4(\hat{\l}''+2\hat{\mu}'')\sqrt{\zeta^2+1}+\kappa_3 (\hat{\l}''+2\hat{\mu}'')\zeta\sqrt{\zeta^2+1}&
\kappa_4\hat{\mu}''\zeta\sqrt{\zeta^2+1}+\kappa_3 \hat{\mu}''\sqrt{\zeta^2+1}
\end{pmatrix}\\[1em]
=&-\frac{h^2}{8}
\begin{pmatrix}
-\dfrac{\hat{\l}''+2\hat{\mu}''}{\zeta^2+1} &-\dfrac{2\hat{\mu}''\zeta}{\zeta^2+1} \\[1em]
\dfrac{2(\hat{\l}''+2\hat{\mu}'')\zeta}{\zeta^2+1}&-\dfrac{\hat{\mu}''}{\zeta^2+1}
\end{pmatrix}
=:T_1.
\end{align*}
\textit{Second term}: the term of order $h^2$ in 
\begin{align*}
&\begin{pmatrix}
\dfrac{1}{\sqrt{\zeta^2+1}} & \dfrac{\zeta}{\sqrt{\zeta^2+1}}\\[1em]
-\dfrac{\zeta}{\sqrt{\zeta^2+1}} & \dfrac{1}{\sqrt{\zeta^2+1}}
\end{pmatrix}\star
\frac{h}{2\ii}
\begin{pmatrix}
\hat{\mu}'\dfrac{\zeta}{\sqrt{\zeta^2+1}}&(\hat{\l}'+2\hat{\mu}')\dfrac{1}{\sqrt{\zeta^2+1}}\\[1em] 
-\hat{\mu}'\dfrac{1}{\sqrt{\zeta^2+1}}&(\hat{\l}'+2\hat{\mu}')\dfrac{\zeta}{\sqrt{\zeta^2+1}}
\end{pmatrix}
\end{align*}
is
\begin{align*}
-\frac{h^2}{4}
&\begin{pmatrix}
\kappa_1 \hat{\mu}''\dfrac{\zeta}{\sqrt{\zeta^2+1}}-\kappa_2 \hat{\mu}''\dfrac{1}{\sqrt{\zeta^2+1}}&
\kappa_1(\hat{\l}''+2\hat{\mu}'')\dfrac{1}{\sqrt{\zeta^2+1}}+\kappa_2 (\hat{\l}''+2\hat{\mu}'')\dfrac{\zeta}{\sqrt{\zeta^2+1}}\\[1em] 
-\kappa_2 \hat{\mu}''\frac{\zeta}{\sqrt{\zeta^2+1}}-\kappa_1  \hat{\mu}''\frac{1}{\sqrt{\zeta^2+1}}&
-\kappa_2(\hat{\l}''+2\hat{\mu}'')\dfrac{1}{\sqrt{\zeta^2+1}}+\kappa_1 (\hat{\l}''+2\hat{\mu}'')\dfrac{\zeta}{\sqrt{\zeta^2+1}}
\end{pmatrix}\\[1em]
=-\frac{h^2}{4}
& \begin{pmatrix}
-\dfrac{\hat{\mu}''}{\zeta^2+1}&0
\\ -0&-\dfrac{\hat{\l}''+2\hat{\mu}''}{\zeta^2+1}
\end{pmatrix}=:T_2.
\end{align*}
\textit{Third term}: the term of order $h^2$ in 
\begin{align*}
&-\begin{pmatrix}
\dfrac{1}{\sqrt{\zeta^2+1}} & \dfrac{\zeta}{\sqrt{\zeta^2+1}}\\[1em]
-\dfrac{\zeta}{\sqrt{\zeta^2+1}} & \dfrac{1}{\sqrt{\zeta^2+1}}
\end{pmatrix}\star
\frac{h^2}{8}
\begin{pmatrix}
\dfrac{\zeta^2(2\hat{\mu}''\zeta^2-\hat{\l}'')-(\hat{\l}''+2\hat{\mu}'')}{(\zeta^2+1)^2\sqrt{\zeta^2+1}}&
\dfrac{\zeta(2\hat{\l}''+5\hat{\mu}'')}{(\zeta^2+1)\sqrt{\zeta^2+1}}\\[2em] 
-\dfrac{\zeta(\hat{\l}''+4\hat{\mu}'')}{(\zeta^2+1)\sqrt{\zeta^2+1}}&
\dfrac{2\zeta^4(\hat{\l}''+2\hat{\mu}'')+\zeta^2(2\hat{\l}''+3\hat{\mu}'')-\hat{\mu}''}{(\zeta^2+1)^2\sqrt{\zeta^2+1}}
\end{pmatrix}
\end{align*}
is
\begin{align*}
&-\frac{h^2}{8(\zeta^2+1)^3}
\begin{pmatrix}
1 & \zeta\\[1em]
-\zeta &1
\end{pmatrix}
\begin{pmatrix}
\zeta^2(2\hat{\mu}''\zeta^2-\hat{\l}'')-(\hat{\l}''+2\hat{\mu}'')&
\zeta(2\hat{\l}''+5\hat{\mu}'')(\zeta^2+1)\\[1em]
-\zeta(\hat{\l}''+4\hat{\mu}'')(\zeta^2+1)&
2\zeta^4(\hat{\l}''+2\hat{\mu}'')+\zeta^2(2\hat{\l}''+3\hat{\mu}'')-\hat{\mu}''
\end{pmatrix}\\[1em]
=&-\frac{h^2}{8(\zeta^2+1)^3}
\begin{pmatrix}
-(\hat{\l}''+2\hat{\mu}'')(\zeta^2+1)^2 & 
2\zeta (\hat{\l}''+2\hat{\mu}'')(\zeta^2+1)^2\\[1em]
-2\zeta\hat{\mu}''(\zeta^2+1)^2 &
-\hat{\mu}''(\zeta^2+1)^2
\end{pmatrix} \\[1em]
=&-\frac{h^2}{8(\zeta^2+1)}
\begin{pmatrix}
-(\hat{\l}''+2\hat{\mu}'') & 
2\zeta (\hat{\l}''+2\hat{\mu}'')\\[1em]
-2\zeta\hat{\mu}'' &
-\hat{\mu}''
\end{pmatrix}=:T_3.
\end{align*}
We also need to take into account the transform of the $h^1$-order
term in $q$ (only to leading order)
\begin{align*}
   Q^{-1}\star q_1\star Q &=\frac{1}{\sqrt{\zeta^2 + 1}}
   \left(\begin{array}{cc} 1 & \zeta
   \\
   -\zeta & 1 \end{array}\right)
   \star\frac{1}{2\ii} \left(\begin{array}{cc}
     0 & \hat{\mu}'-\hat{\lambda}'
   \\
   \hat{\lambda}'-\hat{\mu}' & 0 \end{array}\right)
   \star\frac{1}{\sqrt{\zeta^2+1}}
   \left(\begin{array}{cc} 1 & -\zeta
   \\
   \zeta & 1 \end{array}\right)
\\[1em]
   &= \frac{1}{2\ii} \left(\begin{array}{cc}
         0 & \hat{\mu}'-\hat{\lambda}'
   \\
   \hat{\lambda}'-\hat{\mu}' & 0 \end{array}\right) + {\mathcal O}(h) .
\end{align*}
We require the $h^2$-order terms in $hQ^{-1} \star q_1 \star Q$ in the
further analysis. First, we calculate
\begin{align*}
   h q_1 \star Q &=
   \frac{h}{2\ii}
   \begin{pmatrix}
   0&\hat{\mu}'-\hat{\lambda}'\\[1em]
   \hat{\lambda}'-\hat{\mu}'&0
   \end{pmatrix} \star
   \begin{pmatrix}
   \dfrac{1}{\sqrt{\zeta^2+1}} & -\dfrac{\zeta}{\sqrt{\zeta^2+1}}\\[1em]
   \dfrac{\zeta}{\sqrt{\zeta^2+1}} & \dfrac{1}{\sqrt{\zeta^2+1}}
   \end{pmatrix}
\\[1em]
   &=
   \frac{h}{2\ii}
\begin{pmatrix} 
(\hat{\mu}'-\hat{\lambda}')\star\dfrac{\zeta}{\sqrt{\zeta^2+1}}&
(\hat{\mu}'-\hat{\lambda}')\star\dfrac{1}{\sqrt{\zeta^2+1}}\\[1em]
(\hat{\lambda}'-\hat{\mu}')\star\dfrac{1}{\sqrt{\zeta^2+1}}&
(\hat{\mu}'-\hat{\lambda}')\star\dfrac{\zeta}{\sqrt{\zeta^2+1}}
\end{pmatrix}\\[1em]
   &= \frac{h}{2\ii}
\begin{pmatrix}
(\hat{\mu}'-\hat{\lambda}')\dfrac{\zeta}{\sqrt{\zeta^2+1}}&
(\hat{\mu}'-\hat{\lambda}')\dfrac{1}{\sqrt{\zeta^2+1}}\\[1em]
(\hat{\lambda}'-\hat{\mu}')\dfrac{1}{\sqrt{\zeta^2+1}}&
(\hat{\mu}'-\hat{\lambda}')\dfrac{\zeta}{\sqrt{\zeta^2+1}}
\end{pmatrix}+
\frac{h^2}{4}
\begin{pmatrix}
(\hat{\mu}''-\hat{\lambda}'')\kappa_2&
(\hat{\mu}''-\hat{\lambda}'')\kappa_1\\[1em]
(\hat{\lambda}''-\hat{\mu}'')\kappa_1&
(\hat{\mu}''-\hat{\lambda}'')\kappa_2
\end{pmatrix} =: \mathfrak{T} .
\end{align*}
Then, up to $h^2$-order terms, 
\begin{align*}
   h &Q^{-1} \star q_1 \star Q =
   \begin{pmatrix}
   \dfrac{1}{\sqrt{\zeta^2 + 1}} & \dfrac{\zeta}{\sqrt{\zeta^2 + 1}}
   \\[1em]
   -\dfrac{\zeta}{\sqrt{\zeta^2 + 1}} & \dfrac{1}{\sqrt{\zeta^2 + 1}}
   \end{pmatrix} \star \mathfrak{T}
\\[1em]
   &= \frac{h}{2 \ii}
   \left(\begin{array}{c}
   \dfrac{1}{\sqrt{\zeta^2 + 1}} \star
     (\hat{\mu}' - \hat{\lambda}') \dfrac{\zeta}{\sqrt{\zeta^2 + 1}}
   + \dfrac{\zeta}{\sqrt{\zeta^2 + 1}} \star
       (\hat{\lambda}' - \hat{\mu}') \dfrac{1}{\sqrt{\zeta^2 + 1}} 
   \\[1em]
   \dfrac{1}{\sqrt{\zeta^2 + 1}} \star
       (\hat{\mu}' - \hat{\lambda}') \dfrac{1}{\sqrt{\zeta^2 + 1}}
   + \dfrac{\zeta}{\sqrt{\zeta^2 + 1}} \star
     (\hat{\mu}' - \hat{\lambda}') \dfrac{\zeta}{\sqrt{\zeta^2 + 1}}
   \end{array}\right.
\\[1em]
   &\hspace*{3.5cm} \left.\begin{array}{c}
   \dfrac{1}{\sqrt{\zeta^2 + 1}} \star
       (\hat{\mu}' - \hat{\lambda}') \dfrac{1}{\sqrt{\zeta^2 + 1}}
   + \dfrac{\zeta}{\sqrt{\zeta^2 + 1}} \star
     (\hat{\mu}' - \hat{\lambda}') \dfrac{\zeta}{\sqrt{\zeta^2 + 1}} 
   \\[1em]
   -\dfrac{\zeta}{\sqrt{\zeta^2 + 1}} \star
       (\hat{\mu}' - \hat{\lambda}') \dfrac{1}{\sqrt{\zeta^2 + 1}}
   + \dfrac{1}{\sqrt{\zeta^2 + 1}} \star
     (\hat{\mu}' - \hat{\lambda}') \dfrac{\zeta}{\sqrt{\zeta^2 + 1}}
   \end{array}\right)
\\[1em]
   + &\frac{h^2}{4}
\begin{pmatrix}
\dfrac{1}{\sqrt{\zeta^2+1}}(\hat{\mu}''-\hat{\lambda}'')\kappa_2+ \dfrac{\zeta}{\sqrt{\zeta^2+1}}(\hat{\lambda}''-\hat{\mu}'')\kappa_1 &
\dfrac{1}{\sqrt{\zeta^2+1}}(\hat{\mu}''-\hat{\lambda}'')\kappa_1+ \dfrac{\zeta}{\sqrt{\zeta^2+1}}(\hat{\mu}''-\hat{\lambda}'')\kappa_2 \\[1em]
 -\dfrac{\zeta}{\sqrt{\zeta^2+1}}(\hat{\mu}''-\hat{\lambda}'')\kappa_2+ \dfrac{1}{\sqrt{\zeta^2+1}}(\hat{\lambda}''-\hat{\mu}'')\kappa_1 &
 -\dfrac{\zeta}{\sqrt{\zeta^2+1}}(\hat{\mu}''-\hat{\lambda}'')\kappa_1+ \dfrac{1}{\sqrt{\zeta^2+1}}(\hat{\mu}''-\hat{\lambda}'')\kappa_2
\end{pmatrix} .
\end{align*}
Thus, the $h^2$-order terms in the expression for $h Q^{-1} \star
q_1\star Q$ are
\begin{align*}
&-\frac{h^2}{4}
\begin{pmatrix}
\kappa_1 (\hat{\mu}''-\hat{\lambda}'')\dfrac{\zeta}{\sqrt{\zeta^2+1}}+ \kappa_2 (\hat{\lambda}''-\hat{\mu}'')\dfrac{1}{\sqrt{\zeta^2+1}} &
\kappa_1 (\hat{\mu}''-\hat{\lambda}'')\dfrac{1}{\sqrt{\zeta^2+1}}+ \kappa_2 (\hat{\mu}''-\hat{\lambda}'')\dfrac{\zeta}{\sqrt{\zeta^2+1}} \\[1em]
-\kappa_2 (\hat{\mu}''-\hat{\lambda}'')\dfrac{\zeta}{\sqrt{\zeta^2+1}}+ \kappa_1(\hat{\lambda}''-\hat{\mu}'')\dfrac{1}{\sqrt{\zeta^2+1}} &
-\kappa_2 (\hat{\mu}''-\hat{\lambda}'')\dfrac{1}{\sqrt{\zeta^2+1}}+ \kappa_1 (\hat{\mu}''-\hat{\lambda}'')\dfrac{\zeta}{\sqrt{\zeta^2+1}}
\end{pmatrix}\\[1em]
&+\frac{h^2}{4}
\begin{pmatrix}
\dfrac{1}{\sqrt{\zeta^2+1}}(\hat{\mu}''-\hat{\lambda}'')\kappa_2+ \dfrac{\zeta}{\sqrt{\zeta^2+1}}(\hat{\lambda}''-\hat{\mu}'')\kappa_1 & 
\dfrac{1}{\sqrt{\zeta^2+1}}(\hat{\mu}''-\hat{\lambda}'')\kappa_1+ \dfrac{\zeta}{\sqrt{\zeta^2+1}}(\hat{\mu}''-\hat{\lambda}'')\kappa_2 \\[1em] 
-\dfrac{\zeta}{\sqrt{\zeta^2+1}}(\hat{\mu}''-\hat{\lambda}'')\kappa_2+ \dfrac{1}{\sqrt{\zeta^2+1}}(\hat{\lambda}''-\hat{\mu}'')\kappa_1 & 
-\dfrac{\zeta}{\sqrt{\zeta^2+1}}(\hat{\mu}''-\hat{\lambda}'')\kappa_1+ \dfrac{1}{\sqrt{\zeta^2+1}}(\hat{\mu}''-\hat{\lambda}'')\kappa_2
\end{pmatrix}\\[1em]
=&\frac{h^2}{2}
\begin{pmatrix}
\dfrac{\hat{\mu}''-\hat{\lambda}''}{\zeta^2+1} &0\\ 
0 &\dfrac{\hat{\mu}''-\hat{\lambda}''}{\zeta^2+1}
\end{pmatrix}=:T_4.
\end{align*}
Finally,
\begin{multline*}
   h^2 Q^{-1}\star q_2\star Q
     = h^2 \frac{1}{\sqrt{\zeta^2+1}}
       \left(\begin{array}{cc} 1 & \zeta \\ -\zeta & 1
       \end{array}\right)
    \star \frac14 \left(\begin{array}{cc} \hat{\mu}'' & 0 \\
             0 & \hat{\lambda}''+2\hat{\mu}''\end{array}\right)
    \star \frac{1}{\sqrt{\zeta^2+1}}
    \left(\begin{array}{cc} 1 & -\zeta \\ \zeta & 1 \end{array}
          \right)
\\[1em]
    =\frac{h^2}{4}\frac{1}{\zeta^2+1}
     \left(\begin{array}{cc} (\hat{\lambda}''+2\hat{\mu}'')\zeta^2+\hat{\mu}''&(\hat{\lambda}''+\hat{\mu}'')\zeta \\ (\hat{\lambda}''+\hat{\mu}'')\zeta& \hat{\mu}''\zeta^2+\hat{\lambda}''+2\hat{\mu}''\end{array}\right) =: T_5 .
\end{multline*}
By summing the $h^1$-order terms, we arrive at
\begin{equation}\label{ab}
   Q^{-1}\star q\star Q=\left(\begin{array}{cc}(\hat{\lambda}+2\hat{\mu})(1+\zeta^2)&0\\0&\hat{\mu}(1+\zeta^2)\end{array}\right) +h r,\qq r= 2 \ii\hat{\mu}'\left(\begin{array}{cc}0&-1\\ 1&0\end{array}\right) + {\mathcal O}(h) ,
\end{equation}
where $r = {\mathcal O}(1)$ is the classical zero-order matrix symbol.

Next, we aim to get rid of the off-diagonal terms, $\gamma_1$,
$\beta_1$, while keeping the diagonal terms, $\alpha_1$, $\delta_1$
(which are zero in the Rayleigh case) unchanged. We construct
\[
   B_0 = \ma 0 & b \\ c & 0 \am
\]
such that 
   \begin{multline*}
   (1-hB_0)\star
  \left(\left(\begin{array}{cc}(\hat{\lambda}+2\hat{\mu})(1+\zeta^2)&0\\0&\hat{\mu}(1+\zeta^2)\end{array}\right)
  +
  h \left(\begin{array}{cc}\alpha_1&\beta_1\\\gamma_1&\delta_1\end{array}\right)\right) \star
   (1 + h B_0)
\\[1em] 
   =\left(\begin{array}{cc}(\hat{\lambda}+2\hat{\mu})(1+\zeta^2)+h\alpha_1&0\\0&\hat{\mu}(1+\zeta^2)+h\delta_1\end{array}\right)
  + {\mathcal O}(h^2) .
\end{multline*}
We choose $b,c$ according to
\begin{align*}
(\hat{\lambda}+2\hat{\mu})(1+\zeta^2) b-b\hat{\mu}(1+\zeta^2)=-\beta_1\qq &\Leftrightarrow\qq b=-\frac{\beta_1}{(\hat{\lambda}+\hat{\mu})(1+\zeta^2)} ,\\
\hat{\mu}(1+\zeta^2)c-c(\hat{\lambda}+2\hat{\mu})(1+\zeta^2)=-\gamma_1\qq &\Leftrightarrow\qq c=\frac{\gamma_1}{(\hat{\lambda}+\hat{\mu})(1+\zeta^2)}.
\end{align*}
Hence, using (\ref{ab}), we get
\begin{equation*}
   B_0 = \frac{2 \ii \hat{\mu}'}{(\hat{\lambda}+\hat{\mu})(1+\zeta^2)}
         \ma 0 & 1 \\ 1 & 0 \am .
\end{equation*}

Now, we consider the $h^2$-order terms. Let
$$
   D = \left(\begin{array}{cc}
   (\hat{\lambda} + 2 \hat{\mu})(1 + \zeta^2) & 0 \\
   0 & \hat{\mu}(1 + \zeta^2)\end{array}\right) .
$$
By summing the $h^1$- and $h^2$-order terms, we get
\[
   Q^{-1}\star q\star Q=D +h r_1+h^2r_2 ,\qq
   r_1 = 2 \mathrm{i} \hat{\mu}'
         \left(\begin{array}{cc} 0 & -1 \\ 1 & 0 \end{array}\right) ,
\]
where $r_1 = {\mathcal O}(1)$ is a classical zero-order matrix symbol
and
\begin{align*}
   r_2 &= T_1 + T_2 + T_3 + T_4 + T_5
\\
   &= -\frac{1}{8} \left(\begin{array}{cc}
   -\dfrac{\hat{\l}'' + 2 \hat{\mu}''}{\zeta^2 + 1} &
   -\dfrac{2 \hat{\mu}'' \zeta}{\zeta^2 + 1}
   \\
   \dfrac{2 (\hat{\l}'' + 2 \hat{\mu}'') \zeta}{\zeta^2 + 1} &
   -\dfrac{\hat{\mu}''}{\zeta^2 + 1} \end{array}\right)
   - \frac{1}{4} \left(\begin{array}{cc}
   -\dfrac{\hat{\mu}''}{\zeta^2 + 1} & 0
   \\
   0 & -\dfrac{\hat{\l}'' + 2 \hat{\mu}''}{\zeta^2 + 1}
   \end{array}\right)
\\[1em]
   &\qquad\qquad - \frac{1}{8 (\zeta^2 + 1)} \left(\begin{array}{cc}
   -(\hat{\l}'' + 2 \hat{\mu}'') & 2 \zeta (\hat{\l}'' + 2 \hat{\mu}'')
   \\
   -2 \zeta \hat{\mu}'' & -\hat{\mu}'' \end{array}\right)
\\[1em]
   &\qquad +\frac{1}{2} \left(\begin{array}{cc}
   \dfrac{\hat{\mu}'' - \hat{\lambda}''}{\zeta^2 + 1} & 0
   \\
   0 & \dfrac{\hat{\mu}'' - \hat{\lambda}''}{\zeta^2 + 1}
   \end{array}\right)
   + \frac{1}{4} \frac{1}{\zeta^2 + 1} \left(\begin{array}{cc}
   (\hat{\lambda}'' + 2 \hat{\mu}'') \zeta^2 + \hat{\mu}'' &
   (\hat{\lambda}'' + \hat{\mu}'') \zeta
   \\
   (\hat{\lambda}'' + \hat{\mu}'') \zeta &
   \hat{\mu}'' \zeta^2 + \hat{\lambda}'' + 2 \hat{\mu}''
   \end{array}\right)
\\[1em]
   &=\frac{1}{8} \frac{1}{\zeta^2 + 1} \left(\begin{array}{cc}
   2 (\hat{\lambda}'' + 2 \hat{\mu}'') \zeta^2 - 2 \hat{\l}''
     + 12 \hat{\mu}'' & 0
   \\
   0 & 2 \hat{\mu}'' \zeta^2 + 14 \hat{\mu}'' \end{array} \right) .
\end{align*}
Furthermore,
\begin{multline*}
(1-hB_0)\star \left(D+ h 2 \ii\hat{\mu}'\left(\begin{array}{cc}0&-1\\ 1&0\end{array}\right) + h^2 r_2\right) \star (1 + h B_0)
\\
   = D +h\left(r_1+D\star B_0-B_0\star D\right)+h^2\left(r_2+r_1\star B_0-B_0\star r_1-B_0DB_0\right)+{\mathcal O}(h^3)
   = D + h^2\tilde{r}_2+{\mathcal O}(h^3) ,
\end{multline*}
where 
\begin{multline*}
   h^2\tilde{r}_2=-h^2B_0 \left(\begin{array}{cc}(\hat{\lambda}+2\hat{\mu})(1+\zeta^2)&0\\0&\hat{\mu}(1+\zeta^2)\end{array}\right) B_0
\\
   + h \cdot \mbox{ ``$h$-order term in'' }(D\star B_0-B_0\star D)\\
   + h^2 \cdot \mbox{ ``leading term in'' }(r_1\star B_0-B_0\star r_1)+h^2r_2 .
\end{multline*}

Our goal is to find the diagonal entries of $\tilde{r}_2$. We write
\begin{align*}
   -h^2 B_0 \left(\begin{array}{cc}
    (\hat{\lambda} + 2 \hat{\mu})(1 + \zeta^2) & 0 \\
    0 & \hat{\mu}(1+\zeta^2) \end{array}\right) B_0
   = \frac{4 (\hat{\mu}')^2 h^2}{(\hat{\lambda} + \hat{\mu})^2
           (1 + \zeta^2)}
   \left(\begin{array}{cc}\hat{\mu} & 0 \\
   0 & \hat{\lambda} + 2 \hat{\mu}\end{array}\right) =: T_6
\end{align*}
and
\begin{multline*}
   T_7 := h \cdot\mbox{ ``$h^1$-order term in'' }(D\star B_0-B_0\star D)
\\
   = \frac{h^2}{2\ii} \left(\mstrut{0.55cm}\right.
   \left[(\hat{\l} + 3 \hat{\mu}) (1 + \zeta^2)\right]'_\zeta
   \left[\frac{2\ii\hat{\mu}'}{(\hat{\lambda}+\hat{\mu})(1+\zeta^2)}\right]'_Z-\left[(\hat{\l}+3\hat{\mu})(1+\zeta^2)\right]'_Z\left[\frac{2\ii\hat{\mu}'}{(\hat{\lambda}+\hat{\mu})(1+\zeta^2)}\right]'_\zeta
   \left.\mstrut{0.55cm}\right)
   \left(\begin{array}{cc} 0 & 1 \\ 1 & 0 \end{array}\right) ,
\end{multline*}
which is off-diagonal. Furthermore,
\[
   T_8 := h^2 \cdot \mbox{ ``leading term in'' }(r_1 \star B_0
      - B_0 \star r_1)
   = \frac{4 (\hat{\mu}')^2 h^2}{(\hat{\lambda} + \hat{\mu})(1 + \zeta^2)}
     \left(\begin{array}{cc} 2 & 0 \\ 0 & -2 \end{array}\right) .
\]
It follows that
\[
   T_6 + T_8
   = \frac{4 (\hat{\mu}')^2 h^2}{
             (\hat{\lambda} + \hat{\mu})^2(1 + \zeta^2)}
   \left(\begin{array}{cc} 2 \hl + 3 \hat{\mu} & 0 \\
         0 & -\hat{\lambda} \end{array}\right) .
\]
Finally, we obtain the diagonal terms in $\tilde{r}_2$, that is,
\[
   \tilde{r}_2 = \ma \alpha_2 & 0 \\ 0 & \delta_2 \am
               = r_2 + T_6 + T_7 + T_8 ,
\]
with
\begin{align}
   \alpha_2&= \frac{1}{\zeta^2+1}\left\{\frac{1}{4}\left( (\hat{\lambda}''+2\hat{\mu}'')\zeta^2-\hat{\l}''+6\hat{\mu}''\right)+\frac{4(\hat{\mu}')^2(2\hl+3\hat{\mu})}{(\hat{\lambda}+\hat{\mu})^2} \right\}
\nonumber\\
&=\frac{\hat{\lambda}''+2\hat{\mu}''}{4}+ \frac{1}{\zeta^2+1}\left\{-\frac{1}{2}\hat{\l}''+\hat{\mu}''+\frac{4(\hat{\mu}')^2(2\hl+3\hat{\mu})}{(\hat{\lambda}+\hat{\mu})^2} \right\}
\end{align}
and
\begin{align} \label{eq:delta2}
   \delta_2 = \frac{1}{\zeta^2+1}\left\{\frac{\hat{\mu}''}{4}\left(\zeta^2+7\right)-\frac{4\hl(\hat{\mu}')^2}{(\hat{\lambda}+\hat{\mu})^2} \right\}=\frac{\hm''}{4}+ \frac{1}{\zeta^2+1}\left\{\frac32 \hat{\mu}''-\frac{4\hl(\hat{\mu}')^2}{(\hat{\lambda}+\hat{\mu})^2} \right\} .
\end{align}
If $\tilde{q}$ denotes the previously obtained symbol, then we
construct $B = B_0 + h B_1 + \ldots$, that is, $B_1$ to get rid of the
off-diagonal entries in $\tilde{r}_2$, such that
\begin{multline*}
   \hat{q} \to q^{\rm diag}
   = e^{\ii B/h} \star \hat{q} \star e^{-\ii B/h}
   = \exp{\left(\frac{\ii}{h}{\rm  ad}\,(B)^\star\right)} \hat{q} ,
\\
   q^{\rm diag} = \left(\begin{array}{cc}(\hat{\lambda}+2\hat{\mu})(1+\zeta^2)+h\alpha_1+h^2\alpha_2+\ldots&0\\0&\hat{\mu}(1+\zeta^2)+h\delta_1+h^2\delta_2+\ldots
   \end{array}\right) .
\end{multline*}
The symbol $B_1$ is constructed as $B_0$ before so that diagonal
entries are unchanged. In the above,
\begin{equation} \label{ad}
   \exp{\left(\frac{\ii}{h}{\rm  ad}\,(B)^\star\right)} \hat{q}
   =\exp{\left(\frac{\ii}{h}[B, .]^\star\right)} \hat{q}=\hat{q}+\frac{\ii}{h}[B,\hat{q}]^\star +\frac12 \left(\frac{\ii}{h}\right)^2[B,[B,\hat{q}]^\star]^\star+\ldots
\end{equation}
is a classical symbol, with
\begin{multline*}
\frac{\ii}{h}[B,\hat{q}]^\star=\{B,\hat{q}\}_1-\frac{1}{24}h^2\{B,\hat{q}\}_3+\ldots ,\quad \{B,\hat{q}\}_1=B'_\zeta \hat{q}_Z' -B'_Z\hat{q}'_\zeta ,
\\
\{B,\hat{q}\}_3=B^{(3)}_{\zeta\zeta\zeta} \hat{q}^{(3)}_{ZZZ}-3B^{(3)}_{\zeta\zeta Z} \hat{q}^{(3)}_{ZZ\zeta}+3 B^{(3)}_{\zeta Z Z} \hat{q}^{(3)}_{Z\zeta\zeta}-B^{(3)}_{ZZZ} \hat{q}^{(3)}_{\zeta\zeta\zeta} .
\end{multline*}
{}\hfill \end{proof}

\subsection{Bohr-Sommerfeld quantization rules for multiple wells}
\label{ssec:multiple-wells}

For wells separated from the boundary, the analysis is purely based on
the diagonalized system and, hence, follows the corresponding analysis
for Love waves. That is, we consider operator $H_{0,h,2}$
(cf.~(\ref{diagonal})). We introduce the following assumptions on
$\hat{\mu}$

\medskip

\begin{assumption} \label{assu_boundary well}
There is a $Z^* < 0$ such that $\hat{\mu}'(Z^*) = 0$,
$\hat{\mu}''(Z^*) < 0$ and $\hat{\mu}'(Z) < 0$ for $Z \in
\,]Z^*,0[\,$.
\end{assumption}

\medskip

\begin{assumption} \label{assu2jag}
The function $\hat{\mu}(Z)$ has non-degenerate critical values at a
finite set
$$
   \{ Z_1, Z_2, \cdots, Z_M \}
$$
in $]Z_I,0[$ and all critical points are non-degenerate extrema. None
  of the critical values of $\hat{\mu}(Z)$ are equal, that is,
  $\hat{\mu}(Z_j) \not= \hat{\mu}(Z_k)$ if $j \not= k$.
\end{assumption}

\medskip

We label the critical values of $\hat{\mu}(Z)$ as $E_1 < \ldots < E_M
< \hat{\mu}_I$ and the corresponding critical points by
$Z_1,\cdots,Z_M$. We denote $Z_0 = 0$, $E_0 = \hat{\mu}(Z_0)$ and
$E_{M+1} = \hat{\mu}_I$.

We define a well of order $k$ as a connected component of $\{ Z \in
]Z_I,0[\ :\ \hat{\mu}(Z) < E_k \}$ that is not connected to the
  boundary at $Z = 0$. We refer to the connected component connected
  to the boundary as a half well of order $k$. We denote $J_k =
]E_{k-1},E_k[$, $k = 1,2,3,\cdots$ and let $N_k$ ($\leq k$) be the
  number of wells of order $k$. The set $\{Z \in
  (Z_I,0)\ :\ \hat{\mu}(Z) < E_k\}$ consists of $N_k$ wells and one
  half well
\[
   W_j^k(E) ,\quad j=1,2,\cdots,N_k ,\text{ and }
   \widetilde{W}^k(E) ,\quad
   (\cup_{j=1}^{N_k} W_j^k(E)) \cup \widetilde{W}^k(E)
                                       \subset [Z_I,0[ \, .
\]
The half well $\widetilde{W}^k(E)$ is connected to the boundary at
$Z = 0$.

The semiclassical spectrum mod $o(h^{5/2})$ in $J_k$ is the union of
$N_k+1$ spectra:
\[
   \cup_{j=1}^{N_k} \Sigma_j^k(h) \cup
                    \widetilde{\Sigma}^k(h) .
\]
Here, $\Sigma_j^k(h)$ is the semi-classical spectrum associated to the
well $W_j^k$, and the spectrum $\widetilde{\Sigma}^k(h)$ is the
semiclassical spectrum associated to the half well $\widetilde{W}^k$.

We have Bohr-Sommerfeld rules for separated wells,
\begin{equation} \label{sepa}
   \Sigma_j^k(h) = \{ \mu_\alpha(h)\ :\ E_{k-1} < \mu_\alpha(h)
          < E_k\text{ and }S^{k,j}(\mu_\alpha(h)) = 2\pi h \alpha \} ,
\end{equation}
where $S^{k,j} = S^{k,j}(E) :\ ]E_{k-1},E_k[ \to \mathbb{R}$ admits
the asymptotics in $h$
\[
   S^{k,j}(E) = S^{k,j}_0(E) + h \pi + h^2 S_2^{k,j}(E) + \cdots
\]
and
\begin{equation*}
   \widetilde{\Sigma}^k(h) = \{\nu_\alpha(h)\ :\
   E_{k-1} < \nu_\alpha(h) < E_k\text{ and }
                   \widetilde{S}^k(\nu_\alpha(h)) = 2\pi h \alpha \} ,
\end{equation*}
where $\widetilde{S}^k = \widetilde{S}^k(E) :\ ]E_{k-1},E_k[ \to
\mathbb{R}$ admits the asymptotics in $h$
\[
   \widetilde{S}^k(E) = \frac{1}{2} \tilde{S}^k_0(E)
         + h \widetilde{S}_1^k(E)
              + \frac{1}{2} h^2 \widetilde{S}_2^k(E) + \cdots .
\]

For the explicit forms of $S^{k,j}$ and $\widetilde{S}^k$, we
introduce the classical Hamiltonian
$p_0(Z,\zeta) = \hat{\mu}(Z) (1 + \zeta^2)$ coinciding with the $h^0$
term in $\s^{\rm W}(H_{0,h,2})$. For any $k$, $p_0^{-1}(J_k)$ is a
union of $N_k$ topological annuli $A_j^k$ and a half annulus
$\widetilde{A}^k$. The map $p_0 :\ A_j^k \to J_k$ is a fibration whose
fibers $p_0^{-1}(E) \cap A_j^k$ are topological circles $\g_j^k(E)$
that are periodic trajectories of classical dynamics. The map $p_0
:\ \widetilde{A}^k \to J_k$ is a topological half circle
$\widetilde{\gamma}^k(E)$. If $E \in J_k$ then $p_0^{-1}(E) =
(\cup_{j=1}^{N_k} \g_j^k(E)) \cup \widetilde{\gamma}^k(E)$. The
corresponding classical periods are
\[
   T_j^k(E) = \int_{\g_j^k(E)} |\mathrm{d}t|\quad\text{and}\quad
   \frac12 \widetilde{T}^k(E)
            = \int_{\widetilde{\g}^k(E)} |\mathrm{d}t| .
\]
We let $t$ be the parametrization of $\g^k_j(E)$ by time evolution in
\begin{equation}
\label{traj}
   \frac{\mathrm{d}Z}{\mathrm{d}t} = \partial_\zeta p_0 ,\quad
   \frac{\mathrm{d}\zeta}{\mathrm{d}t} = -\partial_Z p_0
\end{equation}
for a realized energy level $E$.

For a well $W^k_j$ separated from the boundary, we get
\begin{equation} \label{eq:S0}
   S_0^{k,j}(E) = \int_{\g_j^k(E)} \zeta \mathrm{d}Z
\end{equation}
and
\begin{equation} \label{eq:S2}
\begin{split}
   S_2^{k,j}(E) &= -\frac{1}{12} \frac{d}{dE}
        \int_{\g^k_j(E)} \left(E \hat{\mu}''
      - 2 \left(\frac{E}{\hat{\mu}} - 1\right) (\hat{\mu}')^2\right)
                         |\mathrm{d}t|
      - \int_{\g^k_j(E)} \delta_2 |\mathrm{d}t| .
      \end{split}
\end{equation}
Substituting (\ref{eq:delta2}), we obtain
\begin{equation} \label{eq:S2JKL}
   S_2^{k,j}(E) = -\frac{1}{12} \frac{\mathrm{d}}{\mathrm{d}E}J(E)
                  - \frac14 K(E) - L(E) ,
\end{equation}
where
\begin{eqnarray*}
   J(E) &=& 
        \int_{\g^k_j(E)} \left(E \hat{\mu}''
      - 2 \left(\frac{E}{\hat{\mu}} - 1\right) (\hat{\mu}')^2\right)
                         |\mathrm{d}t| ,
\\
   K(E) &=& \int_{\g^k_j(E)} \hat{\mu}'' |\mathrm{d}t| ,
\\
   L(E) &=& \int_{\g^k_j(E)} \frac{\hat{\mu}}{E}
            \left( \frac32 \hat{\mu}''
         - \frac{4 \hl(\hat{\mu}')^2}{(\hat{\lambda} + \hat{\mu})^2}
                            \right) |\mathrm{d}t| .
\end{eqnarray*}
The integrations along the periodic trajectory $\g$ can be changed
into integrations over $]f_-(E),f_+(E)[$, $E \in [E_{k-1},E_k]$, in
the $Z$ coordinate. We get
\begin{equation} \label{formSfirst}
   S_0^{k,j}(E)
     = 2 \int_{f_-(E)}^{f_+(E)}
            \sqrt{\frac{E - \hat{\mu}}{\hat{\mu}}} \mathrm{d}Z
\end{equation}
and
\begin{eqnarray*}
   J(E) &=& \int_{f_-(E)}^{f_+(E)}
     \left( E \hat{\mu}'' - \frac{2 (E - \hat{\mu})}{\hat{\mu}}
            (\hat{\mu}')^2 \right)
                \frac{\mathrm{d}Z}{\sqrt{\hm (E - \hm)}} ,
\\
   K(E) &=& \int_{f_-(E)}^{f_+(E)} \hm''
                \frac{\mathrm{d}Z}{\sqrt{\hm (E - \hm)}} ,
\\
   L(E) &=& \int_{f_-(E)}^{f_+(E)} \frac{\hm}{E}
            \left( \frac32 \hat{\mu}''
        -\frac{4 \hl (\hat{\mu}')^2}{(\hat{\lambda} + \hat{\mu})^2}
                   \right)
                \frac{\mathrm{d}Z}{\sqrt{\hm (E - \hm)}} .
\end{eqnarray*}

For the half well $\widetilde{W}^k$ connected to the boundary, we can
write
\begin{equation} \label{eq:tS0}
   \widetilde{S}^k_0(E)
     = 2 \int_{\widetilde{\gamma}^k(E)} \zeta \mathrm{d}Z
     = 4 \int_{f(E)}^0 \sqrt{\frac{E
                  - \hat{\mu}}{\hat{\mu}}} \mathrm{d}Z ,
\end{equation}
as the integration along the periodic half trajectory $\widetilde{\g}$
can be changed into an integration over $]f(E),0[$, $E \in
[E_{k-1},E_k]$, in the $Z$ coordinate. From (\ref{BSbis}) it follows
that
\begin{equation} \label{eq:tS1}
   \widetilde{S}^k_1(E) = \frac{3\pi}{4}
   + \arctan\left(\frac{4 \hat{\mu}(0)(\hat{\mu}(0) - E)
      \left(1 - \frac{E}{(\hat{\lambda} + 2 \hat{\mu})(0)}
                  \right)^{1/2}}{(2 \hat{\mu}(0) - E)^2
         \left(\frac{E}{\hat{\mu}(0)} - 1\right)^{1/2}}\right) .
\end{equation}
We note that $S_0^{k,j}$ and $\widetilde{S}^k_0$ depend only on
periodic trajectories. Moreover, we note that we only need to consider
the Bohr-Sommerfeld rules for single wells in the analysis of the
inverse problem, because of the fact that the eigenfunctions are
$\mathcal{O}(h^\infty)$ outside the wells.

\section{Unique recovery of $\hat{\mu}$ from the semiclassical
         spectrum}
\label{inverse-R}

Similar as in the case of Love waves, we obtain a trace formula: As
distributions on $J_k$, we have
\begin{align}
   \sum_{\alpha \in \Z} \delta(E - \mu_\alpha(h))
   =& \frac{1}{2\pi h} \sum_{j = 1}^{N_k}
      \sum_{m \in \Z} (-1)^m
      e^{\mathrm{i} m S_0^{k,j}(E) h^{-1}} T^k_j(E)
     (1 + \mathrm{i} m h S_2^{k,j}(E))
\nonumber\\
   &+ \frac{1}{2\pi h} \sum_{m \in \Z}
      e^{\mathrm{i} m \frac12 \widetilde{S}_0^k(E) h^{-1}}
      e^{\mathrm{i} m \widetilde{S}_1^k(E)}
      \left(\widetilde{T}^k(E)
            + h (\widetilde{S}_1^k)'(E)\right)
      \left(1 + \mathrm{i} m h \frac12 \widetilde{S}^k_2(E)\right)
            + o(1)
\label{CdV6boundary}
\end{align}
having replaced $\nu_{\alpha}$ by $\mu_{\alpha}$ in the notation of
the identification of $\widetilde{\Sigma}^k(h)$. We then introduce the
notation
\begin{eqnarray*}
   Z^k_{m,j}(E) &=& \frac{1}{2\pi h} (-1)^m
        e^{\mathrm{i} m S^{k,j}_0(E) h^{-1}} T^k_j(E)
                    (1 + \mathrm{i} m h S^{k,j}_2(E)) ,\quad
   j = 1,\cdots,N_k ,
\\
   Z^k_{m,N_k+1}(E) &=&
        \frac{1}{2\pi h} e^{\mathrm{i} m\widetilde{S}_1^k(E)}
        e^{\mathrm{i} m \frac12 \widetilde{S}^k_0(E) h^{-1}}
   \left(\frac12 \widetilde{T}^k(E)
   + h (\widetilde{S}_1^k)'(E)\right)
     \left(1 + \mathrm{i} m h \frac12 \widetilde{S}^k_2(E)\right) ,
\end{eqnarray*}
for $m \in \Z$. To further unify the notation, we write
\[
   T^k_{N_k+1}(E) := \frac12 \widetilde{T}^k(E) ,\quad
   S^{k,N_k+1}_0(E) := \frac12 \widetilde{S}^k_{0}(E) ,\quad
   S^{k,N_k+1}_1(E) := \widetilde{S}^k_{1}(E) ,\quad
   S^{k,N_k+1}_2(E) := \frac12 \widetilde{S}^k_2(E) .
\]
Then
$$
   Z^k_{m,N_k+1}(E) =
        \frac{1}{2\pi h} e^{\mathrm{i} mS^{k,N_k+1}_{1}(E)}
        e^{\mathrm{i} mS^{k,N_k+1}_{1}(E)h^{-1}}
                     (T^k_{N_k+1}(E) + h (S^{k,N_k+1}_{1})'(E))
          (1 + \mathrm{i} m hS^{k,N_k+1}_{2}(E)) .
$$

\subsection{Separation of clusters}\label{Sep_cluster}

In \cite{dHINZ}, it was proved that there exists a unique eigenvalue
of $H_{0,h}$ below $\hat{\mu}(0)$ for small $h$. This eigenvalue
cannot be related to any well. Therefore, we first separate out this
fundamental mode to continue our presentation. We then follow
\cite[Subsection~5.2]{dHIvdHZ} providing the separation of clusters
for Love modes applying \cite[Lemma~11.1]{CdV2011}. We invoke

\medskip

\begin{assumption} \label{assum5}
For any $k = 1,2,\cdots$ and any $j$ with $1 \leq j < l \leq N_k+1$,
{\em the classical periods (half-period if $j=N_k+1$) $T^k_j(E)$ and
  $T^k_l(E)$ are weakly transverse in $J_k$}, that is, there exists an
integer $N$ such that the $N$th derivative $(T^k_j(E) -
T^k_l(E))^{(N)}$ does not vanish.
\end{assumption}

\medskip

\noindent
As in the case of Love modes, we introduce the sets
\[
   B = \{ E \in J_k\ :\ \exists j \neq l ,\quad T^k_j(E) = T^k_l(E) \} ,
\]
while suppressing $k$ in the notation. By the weak transversality
assumption, it follows that $B$ is a discrete subset of $J_k$.

We let the distributions
\[
   D_h(E) = \sum_{\alpha \in \Z} \delta(E - \mu_\alpha(h))
\]
be given on the interval $J = J_k$ modulo $o(1)$ using
(\ref{CdV6boundary}). Since $J_k \cap
(-\infty,\hat{\mu}(0))=\emptyset$ for any $k$, we can ignore the
lowest eigenvalue $\lambda_0$. These distributions are determined mod
$o(1)$ by the semiclassical spectra mod $o(h^{5/2})$. We denote by
$Z_h$ the finite sum defined by the right-hand side of
(\ref{CdV6boundary}) restricted to $m = 1$,
\[
   Z_h^k(E) = \sum_{j=1}^{N_k+1} Z^k_{1,j}(E) .
\]
Assuming that we already have recovered $\hat{\mu}(0)$, we obtain
$\widetilde{S}_1^k(E)$. By analyzing the microsupport of $D_h$ and
$Z_h$ \cite[Lemmas 12.2 and 12.3]{CdV2011}, we find

\medskip 

\begin{lemma}
Under the weak transversality assumption, the sets $B$ and the
distributions $Z_h^k$ mod $o(1)$ are determined by the distributions
$D_h$ mod $o(1)$.
\end{lemma}

\medskip 

\begin{proof}
As in \cite[Lemma 12.2]{CdV2011}, we do not assume the weak
transversality of the nonprimitive periods $mT_j^k$, $m>1$. For $k=1$,
$Z_h^1(E)$ is associated with only the half well and can be
straightforwardedly recovered.

We now assume that $Z_h^{k-1}(E)$ for $E\in [E_{k-2},E_{k-1}[$ is
already recovered as $Z^{k-1}_{1,N_{k-1}+1}$ (associated with the half
well) has been identified. We write $\tau_1(E) = \inf_j T_j^k(E)$ and
take a maximal interval $K$ with $\inf K=E_{k-1}$ on which $\tau_1$ is
smooth. On $K$, $\tau_1=T_{j_0}^k$ for a unique $j_0$. As in the proof
of \cite[Lemma 12.2]{CdV2011}, we can recover $Z_{1,j_0}^k$ and
$L^k_{1,j_0}$. Then we need to decide whether $j_0$ is equal to $N_k +
1$, which can be done under the weak transversality assumption. If
$j_0 = N_k + 1$, that is, $Z_{1,j_0}^k$ is associated with the half
well, then, with the recovered $\widetilde{S}_1^k(E)$, we can recover
$Z^k_{m,j_0}$ for any $m$. If $j_0 \ne N_{k}+1$, then $Z_{1,j_0}^k$ is
associated with some full well, and $Z^k_{m,j_0}$ for any $m$ can also
be recovered. The proof can be completed following the proof of
\cite[Lemma 12.2]{CdV2011} by continuing this process.
\end{proof}

Similar to \cite[Lemma 12.3]{CdV2011}, we have

\medskip

\begin{lemma}
Assuming that the $S^j$'s are smooth and the $a_j$'s do not vanish,
there is a unique splitting of $Z_h$ as a sum
\[
   Z_h(E) = \frac{1}{2\pi h} \sum_{j=1}^{N_k+1}
         (a_j(E) + h b_j(E)) e^{\mathrm{i} S^j(E)/h} + o(1) .
\]
\end{lemma}

\medskip

\noindent
It follows that the spectrum in $J_k$ mod $o(h^{5/2})$ determines the
actions $S^{k,j}_0(E)$, $S^{k,j}_2(E)$ and $\widetilde{S}^k_0(E)$ and
$\widetilde{S}^k_1(E)$ on $J_k$. This provides the separation of the
data for the $N_k$ wells and the half well. Then, as in
\cite{dHIvdHZ}, we proceed with reconstructing $\hat{\mu}$ from the
functions $S^{k,j}_0(E)$, $S^{k,j}_2(E)$ for any $k$ and $j \leq N_k$
and $\widetilde{S}^k_0(E)$, under

\medskip

\begin{assumption} \label{defect}
The function $\hat{\mu}$ has a generic symmetry defect: If there exist
$X_\pm$ satisfying $\hat{\mu}(X_- ) = \hat{\mu}(X_+) < E$, and for all
$N \in \mn$, $\hat{\mu}^{(N)}(X_-) = (-1)^N \hat{\mu}^{(N)}(X_+)$,
then $\hat{\mu}$ is globally even with respect to $\ha (X_+ + X_-)$ in
the interval $\{ Z\ :\ \hat{\mu}(Z) < E \}$.
\end{assumption}

\begin{figure}[h]
\centering
\includegraphics[width=5 in]{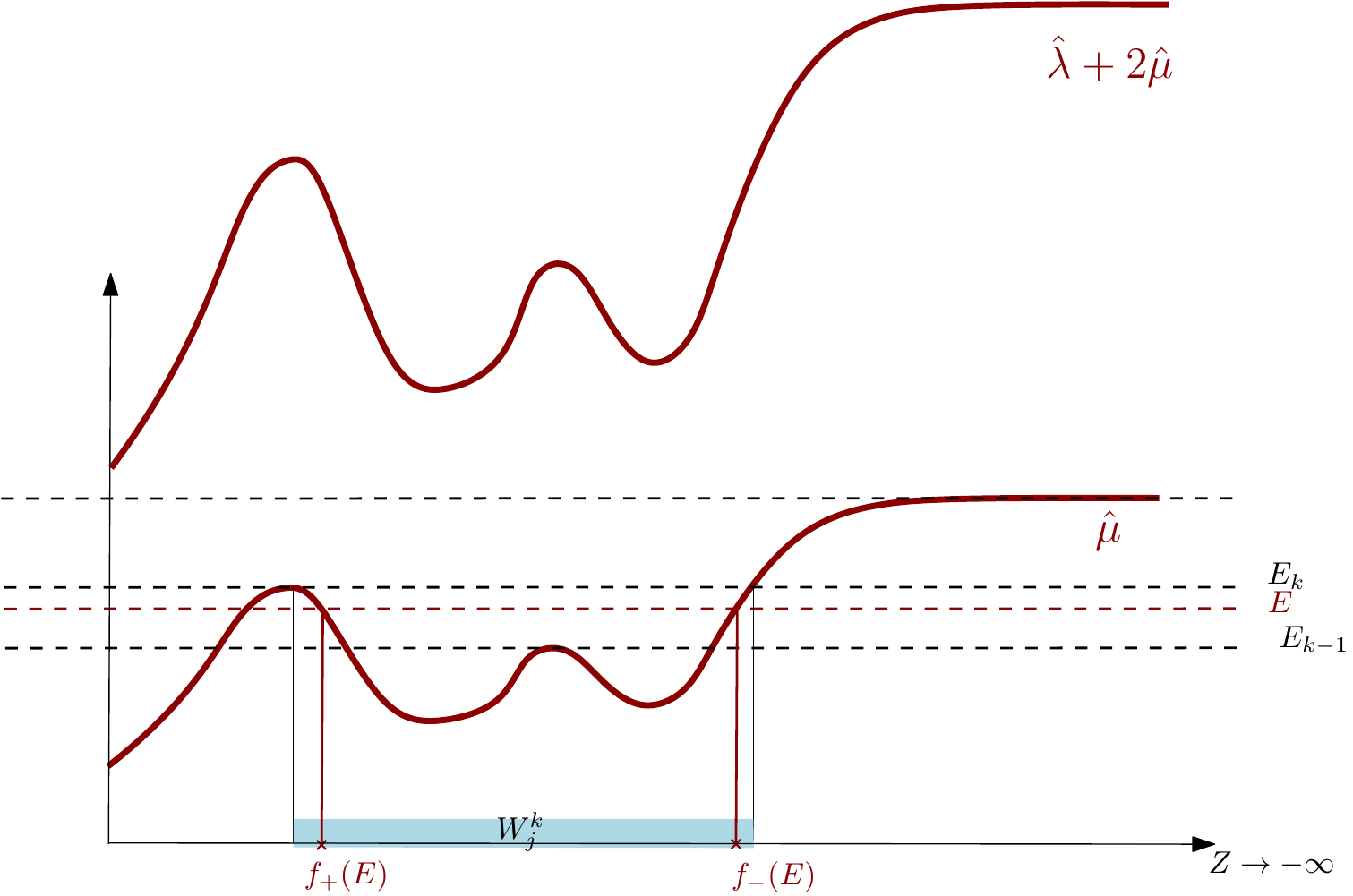}
\caption{Illustration of a well of order $k$ ($N_k = 1$) and
  associated $f_\pm$.}
\label{fpm}
\end{figure}

\subsection{Reconstruction}
\label{ssec:recon}

We note that Assumption~\ref{assu2} is needed here. We summarize the
procedure:

\begin{itemize}
\item
We start by constructing the half well, $\widetilde{W}^1$, that is
connected to the boundary between $E_0$ and $E_1$.
\item
Inductively, we assume that we have already recovered the profile under
$E_{k-1}$. First we reconstruct the half well, $\widetilde{W}^k$, of
order $k$ between $E_{k-1}$ and $E_k$.
\item
We note that $\widetilde{W}^k$ must be a
continuation of the half well $\widetilde{W}^{k-1}$, or be joined with
some well, $W^{k-1}_j$, indexed by $j$ of order $k-1$.
\item
Then we reconstruct a monotonic piece. This can be done as in Section
\ref{S-decreasing} using $\widetilde{S}_0^k$ only.
\item
Secondly, we consider the reconstruction of a full well, $W^k_j$,
separated from the boundary, of order $k$:
\end{itemize}

\medskip

\noindent\textbf{Case I.} The well $W^k_j$ might be a new well. Then
we define the functions $f_\pm :\ [E_{k-1},E_k[ \to I$ so that
    $W^j_k(E) = [f_-(E),f_+(E)]$ for any $E\in [E_{k-1},E_k[$.\\

\noindent\textbf{Case II.} The well $W^k_j$ might also be joining two
wells of order $k-1$, or extending a single well of order $k-1$. Note
that the profile under $E_{k-1}$ has already been recovered. The
smooth joining of two wells can be carried out under
Assumption~\ref{defect}. We consider now functions $f_-(E)$ and
$f_+(E)$ for $E \in [E_{k-1},E_{k}]$ such that $W^k_j$ is the union of
three connected intervals,
$$
   W^k_j(E_k) = [f_-(E_k),f_-(E_{k-1})[ \, \cup \,
   [f_-(E_{k-1}),f_+(E_{k-1})] \, \cup \, ]f_+(E_{k-1}),f_+(E_{k})] .
$$
For an illustration, see Figure~\ref{fpm}.

\medskip

\noindent
For either case, we define
\[
   \Phi(E) = f_+'(E) - f_-'(E) ,\qquad
   \Psi(E) = \frac{1}{f_+'(E)} - \frac{1}{f_-'(E)}.
\]
The recovery goes through explicit reconstruction of the entire
profile following from the gluing procedure as outlined in
\cite[Section~5.4]{dHIvdHZ}. As in the case of the Love modes, the
function $\Phi$ can be recovered from $S_0^{k,j}(E)$, on
$]E_{k-1},E_k[$. From $S^{k,j}_2(E)$, we recover
\begin{multline} \label{operatorA}
   \mathcal{B} \Psi(E) = \int_{E_{k-1}}^E
   \left( (7 E - 6 u) \Psi'(u)
       - 2 \left(\frac{E}{u} - 1\right) \Psi(u) \right)
                        \frac{\mathrm{d}u}{\sqrt{u (E - u)}}
\\
   - \int_{E_{k-1}}^E \left( 36 \Psi'(u)
             - 24 \sigma \frac{1}{u} \Psi(u) \right)
          \arctan\sqrt{\frac{{E - u}}{u}} \mathrm{d}u ,
\end{multline}
$E_{k-1} < E < E_k$, where $\sigma = 8 \nu (1 - 2 \nu)$. This is
established in Appendix~\ref{ssec:A1}. We introduce operator $T$
according to
\[
   T g(E) = \int_{E_{k-1}}^E\frac{g(u)}{\sqrt{E-u}} \mathrm{d}u .
\]
In Appendix~\ref{ssec:A2}, upon setting $E = z^2$, we prove that
\begin{equation} \label{equation_TA}
   \frac{2}{\pi} z^2 \frac{\mathrm{d}^3}{\mathrm{d}z^3}
           (T \circ \mathcal{B} \Psi)(z^2)
   = 16 z^6 \Psi'''(z^2) - 192 z^4 \Psi''(z^2)
        + 96 (2 - \sigma) z^2 \Psi'(z^2) - 96 \sigma \Psi(z^2) .
\end{equation}
That is, we end up with a third-order inhomogeneous ordinary
differential equation for $\Psi(z^2)$ nonsingular on the interval
$[\sqrt{E_{k-1}},\sqrt{E_k}[ \,$. This equation needs to be
    supplemented with ``initial'' conditions:

\medskip

\noindent
For \textbf{Case I}, $\Psi(E_{k-1})$ and the asymptotic behaviors of
$\Psi'(E)$ and $\Psi''(E)$ for $E$ in a neighborhood of $E_{k-1}$ can
be extracted from $T \circ \mathcal{B} \Psi(E)$ and its derivatives at
$E_{k-1}$. Clearly, $\Psi(E_{k-1}) = 0$. Using the derivatives
evaluated in Appendix~\ref{ssec:A2} and
\[
   \Psi(E_{k-1}) = 0 ,\quad
   \lim_{E \downarrow E_{k-1}} \sqrt{E - E_{k-1}} \Psi'(E)
                           = \sqrt{2 \hat{\mu}''(Z_{k-1})} ,
\]
we obtain, for $E>E_{k-1}$ close to $E_{k-1}$,
\[
   \lim_{E \downarrow E_{k-1}} \left(4  E \Psi'(E)
                - \frac{2}{\pi}
   \frac{\mathrm{d}}{\mathrm{d}z}
                (T \circ \mathcal{B} \Psi)(z^2)\right) = 0
\]
yielding the asymptotic behavior of $\Psi'(E)$, and
\[
   \lim_{E \downarrow E_{k-1}} \left(-108 E^{1/2} \Psi'(E)
                + 8 E^{3/2} \Psi''(E) 
        - \frac{\mathrm{d}^2}{\mathrm{d}z^2}
                (T \circ \mathcal{B} \Psi)(z^2)\right) = 0
\]
yielding the asymptotic behavior of $\Psi''(E)$. With these, the
solution to the third-order inhomogeneous ordinary differential
equation is unique.

\medskip

\noindent
For \textbf{Case II}, $\Psi(E_{k-1})$, $\Psi'(E_{k-1})$ and
$\Psi''(E_{k-1})$ are all nonsingular. That is, if $E_{k-1}$ is a
local maximum, $\Psi$ and all its derivatives are smooth from above
and below, and therefore they can be recovered from the reconstruction
on $J_{k-1}$ through one-sided limits. We note that in case $E_{k-1}$
is a local maximum in the middle of two wells in $J_{k-1}$ the two
different $\Psi$s for each well are not smooth below $E_{k-1}$, but it
does not matter as in $J_k$ (above $E_{k-1}$) we use $f_\pm$ from the
monotonically increasing slopes continued from $J_{k-1}$. Thus the
solution to the third-order inhomogeneous ordinary differential
equation is also unique.

\medskip

\noindent
With the recovery of $\Phi$ and $\Psi$ we can recover $f_{\pm}$ and
then $\hat{\mu}$ as in the case of Love modes, again, subject to a
gluing procedure.

\section*{Acknowledgement}

The authors thank Y. Colin de Verdi\`ere for invaluable
discussions. M.V.d.H. gratefully acknowledges support from the Simons
Foundation under the MATH $+$ X program and from NSF under grant
DMS-1815143.

\appendix

\section{Recovery of $\Psi$}

\subsection{Proof of \eqref{operatorA}}
\label{ssec:A1}

We start with the expressions for $J(E)$, $K(E)$ and $L(E)$ in
Subsection~\ref{ssec:multiple-wells}. We apply a change of variable of
integration and obtain,
\begin{align*}
   J(E) =& \int_{E_{k-1}}^E \left(
   E \frac{\mathrm{d}}{\mathrm{d}u}
     \left(\frac{1}{f_+'(u)} - \frac{1}{f_-'(u)}\right)
     - 2 \left(\frac{E}{u} - 1\right)
       \left(\frac{1}{f_+'(u)} - \frac{1}{f_-'(u)}\right)\right)
                   \frac{\mathrm{d}u}{\sqrt{u (E - u)}} + J_{k-1}(E) ,
\\
   K(E) =& \int_{E_{k-1}}^E \frac{\mathrm{d}}{\mathrm{d}u}
           \left(\frac{1}{f_+'(u)} - \frac{1}{f_-'(u)}\right)
                   \frac{\mathrm{d}u}{\sqrt{u (E - u)}} + K_{k-1}(E) ,
\\
   L(E) =& \int_{E_{k-1}}^E \frac{u}{E} \left(
     \frac32 \frac{\mathrm{d}}{\mathrm{d}u}
       \left(\frac{1}{f_+'(u)} - \frac{1}{f_-'(u)}\right)
   - 4 \left(\frac{\hl(f_+(u))}{(\hat{\lambda}(f_+(u)) + u)^2}
       \frac{1}{f_+'(u)}
   - \frac{\hl(f_-(u))}{(\hat{\lambda}(f_-(u)) + u)^2}
       \frac{1}{f_-'(u)}\right) \right)
\\
   &\qquad\quad    \frac{\mathrm{d}u}{\sqrt{u (E - u)}} + L_{k-1}(E) .
\end{align*}
For $\textbf{Case I}$ (cf. Subsection~\ref{ssec:recon}), $J_{k-1}(E)$,
$K_{k-1}(E)$ and $L_{k-1}(E)$ vanish. For \textbf{Case II},
$J_{k-1}(E)$, $K_{k-1}(E)$ and $L_{k-1}(E)$ are related to the profile
on $[f_-(E_{k-1}),f_+(E_{k-1})]$:
\begin{eqnarray}
  J_{k-1}(E) &=& \int_{Z_-}^{Z_+}
  \left(E \hat{\mu}''(Z) - 2 \left(\frac{E}{\hat{\mu}(Z)} - 1\right)
        (\hat{\mu}'(Z))^2\right)
               \frac{\mathrm{d}Z}{\sqrt{\hm(Z)(E - \hm(Z))}} ,
\\
   K_{k-1}(E) &=& \int_{Z_-}^{Z_+} \hm''
               \frac{\mathrm{d}Z}{\sqrt{\hm(Z)(E - \hm(Z))}} ,
\\
   L_{k-1}(E) &=& \int_{Z_-}^{Z_+} \frac{\hm}{E}
       \left(\frac32 \hat{\mu}''
   - \frac{4 \hl(\hat{\mu}')^2}{(\hat{\lambda} + \hat{\mu})^2}
               \right) \frac{\mathrm{d}Z}{\sqrt{\hm(Z)(E - \hm(Z))}} ,
\end{eqnarray}
where $Z_- = f_-(E_{k-1})$ and $Z_+ = f_+(E_{k-1})$. These are already
known.

We find that
\begin{align*}
   K(E) - K_{k-1}(E) =&\ 2 \frac{\mathrm{d}}{\mathrm{d}E}
       \int_{E_{k-1}}^E (E - u) \Psi'(u)
                           \frac{\mathrm{d}u}{\sqrt{u (E - u)}} ,
\\
   L(E) - L_{k-1}(E) =&\ 3 \frac{\mathrm{d}}{\mathrm{d}E}
       \int_{E_{k-1}}^E \left( \Psi'(u)
                - 2 \sigma \frac{1}{u} \Psi(u) \right)
                    \arctan\sqrt{\frac{{E - u}}{u}} \mathrm{d}u .
\end{align*}
Following \cite[Lemma 13.1]{CdV2011}, we introduce an operator
$\mathcal{B}$ defined by
\begin{multline}
   \mathcal{B} \Psi(E) = \int_{E_{k-1}}^E
   \left( (7 E - 6 u) \Psi'(u)
       - 2 \left(\frac{E}{u} - 1\right) \Psi(u) \right)
                        \frac{\mathrm{d}u}{\sqrt{u (E - u)}}
\\
   - \int_{E_{k-1}}^E \left( 36 \Psi'(u)
             - 24 \sigma \frac{1}{u} \Psi(u) \right)
          \arctan\sqrt{\frac{{E - u}}{u}} \mathrm{d}u .
\end{multline}
Using (\ref{eq:S2JKL}), we have established that the derivative of
$\mathcal{B} \Psi$ can be recovered from $S_2^{k,j}$. Then
$\mathcal{B} \Psi$ itself can be recovered using $\mathcal{B}
\Psi(E_{k-1}) = \pi \sqrt{2 \hat{\mu}''(Z_{k-1}) E_{k-1}}$.

\subsection{Proof of \eqref{equation_TA}}
\label{ssec:A2}

We have
\begin{equation}
   (T \circ \mathcal{B} \Psi)(E) = I_1(E) + I_2(E)
                       + I_3(E) + \sigma I_4(E) ,
\end{equation}
where
\begin{eqnarray}
   I_1(E) &=& \int_{E_{k-1}}^E \frac{u}{\sqrt{E - u}}
   \int_{E_{k-1}}^u \left(7 \Psi'(v) - \frac{2}{v} \Psi(v)\right)
   \frac{1}{\sqrt{v}} \frac{1}{\sqrt{u - v}} \,
                       \mathrm{d}v \mathrm{d}u ,
\\
   I_2(E) &=& \int_{E_{k-1}}^E \frac{1}{\sqrt{E - u}}
   \int_{E_{k-1}}^u \left(-6 v \Psi'(v) + 2 \Psi(v)\right)
   \frac{1}{\sqrt{v}} \frac{1}{\sqrt{u - v}} \,
                       \mathrm{d}v \mathrm{d}u ,
\\
   I_3(E) &=& -36 \int_{E_{k-1}}^E \frac{1}{\sqrt{E - u}}
   \int_{E_{k-1}}^u \arctan\sqrt{\frac{{u - v}}{v}} \Psi'(v) \,
                       \mathrm{d}v \mathrm{d}u ,
\\
   I_4(E) &=& 24 \int_{E_{k-1}}^E \frac{1}{\sqrt{E - u}}
   \int_{E_{k-1}}^u \frac{1}{v} \arctan\sqrt{\frac{{u - v}}{v}}
            \Psi(v) \, \mathrm{d}v \mathrm{d}u .
\end{eqnarray}
Upon integration by parts, we obtain
\begin{align*}
   I_3(E) =& -36 \int_{E_{k-1}}^E \frac{1}{u} \sqrt{E - u}
   \left(\int_{E_{k-1}}^u \sqrt{v} \Psi'(v)
        \frac{\mathrm{d}v}{\sqrt{u - v}}\right) \mathrm{d}u
\\
   =& -36 E \int_{E_{k-1}}^E \frac{1}{u} \frac{1}{\sqrt{E - u}}
   \left(\int_{E_{k-1}}^u \sqrt{v} \Psi'(v)
        \frac{\mathrm{d}v}{\sqrt{u - v}}\right) \mathrm{d}u
\\
   &\qquad\qquad
      + 36 \int_{E_{k-1}}^E \frac{1}{\sqrt{E - u}} \left(
               \int_{E_{k-1}}^u \sqrt{v} \Psi'(v)
        \frac{\mathrm{d}v}{\sqrt{u - v}}\right) \mathrm{d}u ,
\\
   I_4(E) =& 24 \int_{E_{k-1}}^E \frac{1}{u} \sqrt{E-u}
       \left(\int_{E_{k-1}}^u \frac{1}{\sqrt{v}} \Psi(v)
        \frac{\mathrm{d}v}{\sqrt{u - v}}\right) \mathrm{d}u
\\
   =& 24 E \int_{E_{k-1}}^E \frac{1}{u} \frac{1}{\sqrt{E - u}}
   \left(\int_{E_{k-1}}^u \frac{1}{\sqrt{v}} \Psi(v)
        \frac{\mathrm{d}v}{\sqrt{u - v}}\right) \mathrm{d}u
   - 24 \int_{E_{k-1}}^E \frac{1}{\sqrt{E - u}}
       \left(\int_{E_{k-1}}^u \frac{1}{\sqrt{v}} \Psi(v)
        \frac{\mathrm{d}v}{\sqrt{u - v}}\right) \mathrm{d}u .
\end{align*}
Now, we use some calculus
\begin{eqnarray*}
   \int_{E_{k-1}}^E \frac{1}{\sqrt{E - u}}
   \left(\int_{E_{k-1}}^u g(v) \frac{\mathrm{d}v}{\sqrt{u - v}}\right)
           \mathrm{d}u &=& \pi \int_{E_{k-1}}^E g(v) \mathrm{d}v ,
\\
   \int_{E_{k-1}}^E \frac{u}{\sqrt{E - u}}
   \left(\int_{E_{k-1}}^u g(v) \frac{\mathrm{d}v}{\sqrt{u - v}}\right)
           \mathrm{d}u &=& \frac{\pi}{2}
                       \int_{E_{k-1}}^E (v + E) g(v) \mathrm{d}v ,
\\
   \int_{E_{k-1}}^E \frac{1}{u} \frac{1}{\sqrt{E - u}} \left(
   \int_{E_{k-1}}^u g(v) \frac{\mathrm{d}v}{\sqrt{u - v}}\right)
           \mathrm{d}u &=& \pi
          \int_{E_{k-1}}^E \frac{1}{\sqrt{E v}} g(v) \mathrm{d}v
\end{eqnarray*}
and get
\begin{eqnarray*}
   I_1(E) &=& \frac{\pi}{2} \int_{E_{k-1}}^E (v + E) \left(7 \Psi'(v)
        - \frac{2}{v} \Psi(v)\right) \frac{1}{\sqrt{v}} \mathrm{d}v ,
\\
   I_2(E) &=& \pi \int_{E_{k-1}}^E \left(-6 v \Psi'(v)
                   + 2\Psi(v)\right) \frac{1}{\sqrt{v}} \mathrm{d}v ,
\\
   I_3(E) &=& 36 \pi \int_{E_{k-1}}^E \left(v - \sqrt{E v}\right)
                            \Psi'(v) \frac{1}{\sqrt{v}} \mathrm{d}v ,
\\
   I_4(E) &=& 24 \pi \int_{E_{k-1}}^E \left(
   \frac{\sqrt{E}}{\sqrt{v}} - 1\right) \Psi(v)
                                     \frac{1}{\sqrt{v}} \mathrm{d}v .
\end{eqnarray*}
We insert $E=z^2$, when trivially
\begin{eqnarray*}
   \frac{2}{\pi} I_1(z^2) &=& \int_{E_{k-1}}^{z^2}
        (v + z^2) \left(7 \Psi'(v) - \frac{2}{v} \Psi(v)\right)
                                     \frac{1}{\sqrt{v}} \mathrm{d}v ,
\\
   \frac{2}{\pi} I_2(z^2) &=& 2 \int_{E_{k-1}}^{z^2}
        \left(-6 v \Psi'(v) + 2 \Psi(v)\right)
                                     \frac{1}{\sqrt{v}} \mathrm{d}v ,
\\
   \frac{2}{\pi} I_3(z^2) &=& 72 \int_{E_{k-1}}^{z^2}
                     \left(\sqrt{v} - z\right) \Psi'(v) \mathrm{d}v ,
\\
   \frac{2}{\pi}I_4(z^2) &=& 48 \int_{E_{k-1}}^{z^2} \left(
        \frac{z}{v} - \frac{1}{\sqrt{v}}\right) \Psi(v) \mathrm{d}v .
\end{eqnarray*}
By tedious calculations, we then find that
\begin{eqnarray*}
   \frac{2}{\pi} \frac{\mathrm{d}}{\mathrm{d}z} I_1(z^2) &=&
                 28 z^2 \Psi'(z^2) - 8 \Psi(z^2)
   + \int_{E_{k-1}}^{z^2} 2 z \left(7 \Psi'(v)
        - \frac{2}{v} \Psi(v)\right) \frac{1}{\sqrt{v}} \mathrm{d}v ,
\\
   \frac{2}{\pi} \frac{\mathrm{d}^2}{\mathrm{d}z^2} I_1(z^2) &=&
   68 z \Psi'(z^2) + 56 z^3 \Psi''(z^2) -8 \frac{1}{z} \Psi(z^2)
   + \int_{E_{k-1}}^{z^2} 2\left(7 \Psi'(v)
        - \frac{2}{v} \Psi(v)\right) \frac{1}{\sqrt{v}} \mathrm{d}t ,
\\
   \frac{2}{\pi} \frac{\mathrm{d}^3}{\mathrm{d}z^3} I_1(z^2) &=&
         112 z^4 \Psi'''(z^2) + 304 z^2 \Psi''(z^2) + 80 \Psi'(z^2)
\end{eqnarray*}
and
\begin{eqnarray*}
   \frac{2}{\pi} \frac{\mathrm{d}}{\mathrm{d}z} I_2(z^2) &=&
                                   -24 z^2 \Psi'(z^2) + 8 \Psi(z^2) ,
\\
   \frac{2}{\pi} \frac{\mathrm{d}^2}{\mathrm{d}z^2} I_2(z^2) &=&
                              -48 z^3 \Psi''(z^2) - 32 z \Psi'(z^2) ,
\\
   \frac{2}{\pi} \frac{\mathrm{d}^3}{\mathrm{d}z^3} I_2(z^2) &=&
         -96 z^4 \Psi'''(z^2) - 208 z^2 \Psi''(z^2) - 32 \Psi'(z^2)
\end{eqnarray*}
and
\begin{eqnarray*}
   \frac{2}{\pi} \frac{\mathrm{d}}{\mathrm{d}z} I_3(z^2) &=&
                      -72 \int_{E_{k-1}}^{z^2} \Psi'(v) \mathrm{d}v ,
\\
   \frac{2}{\pi} \frac{\mathrm{d}^2}{\mathrm{d}z^2} I_3(z^2) &=&
                                                  -144 z \Psi'(z^2) ,
\\
   \frac{2}{\pi} \frac{\mathrm{d}^3}{\mathrm{d}z^3} I_3(z^2) &=&
                              -288 z^2 \Psi''(z^2) - 144 \Psi'(z^2)
\end{eqnarray*}
and
\begin{eqnarray*}
   \frac{2}{\pi} \frac{\mathrm{d}}{\mathrm{d}z} I_4(z^2) &=&
            48 \int_{E_{k-1}}^{z^2} \frac{1}{v} \Psi(v) \mathrm{d}v ,
\\
   \frac{2}{\pi}\frac{\mathrm{d}^2}{\mathrm{d}z^2} I_4(z^2) &=&
                                           96 \frac{1}{z} \Psi(z^2) ,
\\
   \frac{2}{\pi}\frac{\mathrm{d}^3}{\mathrm{d}z^3}I_4(z^2) &=&
                        192 \Psi'(z^2) - 96 \frac{1}{z^2} \Psi(z^2) .
\end{eqnarray*}
These identities lead us to \eqref{equation_TA}.

\bibliographystyle{siam}

\begin{thebibliography}{1}

\bibitem{AbramowitzStegun1965}
M. Abramowitz and I.A. Stegun.
\newblock {\em Handbook of Mathematical Functions}.
\newblock Dover, New York, 2002.

\bibitem{Brocher2005} M. Brocher, {\em Empirical Relations between
  Elastic Wavespeeds and Density in the Earth's Crust},
  Bull. Seism. Soc. Am. \textbf{95}, (2005), pp.2081--2092.

\bibitem{CdV2005}
Y.~Colin de~{V}erdi\`ere.
\newblock Bohr-{S}ommerfeld {R}ules to {A}ll {O}rders.
\newblock {\em Ann. Henri Poincar{\'e}}, 6:925--936, 2005.

\bibitem{CdV2011}
Y.~Colin de~{V}erdi\`ere.
\newblock A semi-classical inverse problem {II}: reconstruction of the
  potential.
\newblock In {\em Geometric aspects of analysis and mechanicss}, volume 292,
  pages 97--119. Progr. Math., 2011.

\bibitem{dHINZ} M.V. de Hoop, A. Iantchenko, G. Nakamura, J. Zhai,
  {\em Semiclassical analysis of elastic surface waves},
  arXiv:1709.06521.

\bibitem{dHIvdHZ} M.V. de Hoop, A. Iantchenko, R. D. Van der Hilst,
  J. Zhai, {\em Semiclassical inverse spectral problem for elastic
    Love surface waves in isotropic media}, arXiv:1908.10529.

\bibitem{Campillo1} A. Derode, E. Larose, M. Tanter, J. de Rosny,
  A. Tourin, M. Campillo and M. Fink, {\em Recovering the Green's
    function form field-field correlations in an open scattering
    medium}, J. Acoust. Soc. Am. \textbf{113}, (2004), pp. 2973--2976.

\bibitem{DimassiSjostrand1999} M.~Dimassi and J.~Sj{\"o}strand,
  Spectral asymptotics in the semiclassical limit, {\em London
    Mathematical Society Lecture Notes Series 268}, Cambridge
  University Press (1991).

\bibitem{HelfferSjostrand1990} B.~Helffer and J.~Sj{\"o}strand,
  Analyse semi-classique de l'{\'e}quation de {H}arper, {II},
  {C}omportement semi-classique pr{\`e}s d'un rationnel, {\em
    M{\'e}m. Soc. Math. France}, 40, 1990.

\bibitem{Martinez2002}
A.~Martinez.
\newblock {\em An {I}ntroduction to {S}emiclassical and {M}icrolocal
  {A}nalysis}.
\newblock Springer, Berlin, 2002.

\bibitem{Campillo2} N. Shapiro and M. Campillo, {\em Emergence of
  broadband Rayleigh waves from correlations of the ambient seismic
  noise}, Geophys. Res. Lett. \textbf{31}, (2004), L07614.

\bibitem{Campillo3} N. Shapiro, M. Campillo, L. Stehly and
  M. Ritzwoller, {\em High resolution surface wave tomography from
    ambient seismic noise}, Science \textbf{307}, (2005), 1615.

\bibitem{KennettWoodhouse1978} B.L.N. Kennett and J.H. Woodhouse, {\em
  On high-frequency spheroidal modes and structure of the upper
  mantle}, Geophys. J.R. astr. Soc. \textbf{55}, (1978), pp.333--350.

\bibitem{Woodhouse1978} J.H. Woodhouse, {\em Asymptotic results for
  elastodynamic propagator matrices in plane-stratified and
  spherically-stratified earth models},
  Geophys. J.R. astr. Soc. \textbf{54}, (1978), pp.263--280.

\end{thebibliography}

\end{document}